\documentclass{siamart171218}
\usepackage[utf8]{inputenc}
\usepackage{amssymb}
\usepackage{enumitem}
\usepackage{amssymb}
\usepackage{mathrsfs}
\usepackage{amsfonts}
\usepackage{graphicx}
\usepackage{subfigure}
\usepackage{xcolor}
 \usepackage{epstopdf}
 \usepackage{boondox-cal}
\usepackage{txfonts}
\usepackage{dsfont}
\usepackage{txfonts}
\usepackage{amssymb}
\usepackage{amsmath}
 \usepackage{epstopdf}
\usepackage{amsmath,amssymb}
\usepackage{mathrsfs,color}
\numberwithin{theorem}{section}
\newsiamremark{remark}{Remark}
\newsiamremark{example}{Example}

\def\theequation{\arabic{section}.\arabic{equation}}
\numberwithin{equation}{section}

\newcommand{\vecx}{\boldsymbol{x}}

\newcommand{\TheTitle}{Implicit integration factor method coupled with Pad\'{e} approximation strategy
 for nonlocal Allen-Cahn equation}

\title{{\TheTitle}\thanks{This work was partially supported by the
 National Natural Science Foundation of China
 (Grant Nos. 12461069, 11961057), the Science and
Technology Project of Guangxi (Grant No. GuikeAD21220114).}}
\author{
Yuxin Zhang\thanks{School of Mathematics and Statistics,
Guangxi Normal University, Guilin 541006, China. (E-mail: zhangyuxin2006@163.com).}
\and
Hengfei Ding\thanks{1.School of Mathematics and Statistics, Guangxi Normal University, Guilin 541006, China;
2.The Center for Applied
Mathematics of Guangxi (GXNU), Guilin 541006, China;
3.Guangxi Colleges and Universities Key Laboratory of Mathematical Model and Application (GXNU),
 Guilin 514006, China. (E-mail:dinghf05@163.com).}}

\begin{document}
\maketitle
\begin{abstract}
 The space nonlocal Allen-Cahn equation is a famous example of fractional reaction-diffusion equations. It is also an extension of the classical Allen-Cahn equation, which is widely used in physics to describe the phenomenon of two-phase fluid flows.
 Due to the nonlocality of the nonlocal operator, numerical solutions to
 these equations face considerable challenges.
 It is worth noting that whether we use low-order or
 high-order numerical differential
  formulas to approximate the operator, the corresponding matrix is always
  dense, which implies that the storage space and computational cost
  required for the former and the latter are the same.
 However, the higher-order formula can significantly
  improve the accuracy of the numerical scheme.
   Therefore, the primary goal of this paper is to construct a high-order numerical formula that approximates the nonlocal operator.
  To reduce the time step limitation in existing numerical algorithms, we employed a technique combining the compact integration factor method with the Pad\'{e} approximation strategy to discretize the time derivative.
A novel high-order numerical scheme, which satisfies both the maximum principle and energy stability for the space nonlocal Allen-Cahn equation, is proposed.
Furthermore, we provide a detailed error analysis of the differential scheme, which shows that its convergence order is $\mathcal{O}\left(\tau^2+h^6\right)$.
 Especially, it is worth mentioning that the fully implicit scheme with sixth-order accuracy
 in spatial has never been proven to maintain the maximum principle and energy stability before.
Finally, some numerical experiments are carried out to demonstrate
 the efficiency of the proposed method.
\end{abstract}


\section{Introduction}\label{sec: introduction}
In this paper, we focus on numerical analysis of the phase field models involving space nonlocaity.
More specifically,
we consider the following space fractional Allen-Cahn equation
\begin{equation}\label{eq.1.1}
\begin{aligned}\displaystyle
\partial_tu-\epsilon^2\mathscr{L}^\gamma u
=-F'(u):=f(u),\;\;(\vecx,t)\in\Omega\times\left(0,T\right],
\end{aligned}
\end{equation}
with the initial condition
\begin{equation*}
\begin{aligned}\displaystyle
u(\vecx,0)=u_0(\vecx),\;\vecx\in\overline{\Omega},
\end{aligned}
\end{equation*}
and the homogeneous Dirichlet boundary condition
\begin{equation*}
\begin{aligned}\displaystyle
u(\vecx,t)|_{\partial\Omega}=0,\;t\in(0,T],
\end{aligned}
\end{equation*}
where the function $F(u)$ is bistable, e.g.,
$F(u)=\frac{1}{4}\left(u^2-1\right)^2$ represents a double-well potential. Here
$\Omega\subset \mathds{R}^d \; (d = 1, 2, 3)$ is an open and bounded Lipschitz domain,
$u : \overline{\Omega} \times [0,\infty) \rightarrow \mathds{R}$ is the unknown function,
$\varepsilon > 0$ is the interfacial parameter representing
the width of the transition regions \cite{Nec2008}.
The operator $\mathcal{L}_\gamma$ denotes the Riesz-type fractional derivative of order
$\gamma\in(1,2]$ in space, which is a typical example of nonlocal operators and
the one-dimensional case is defined by \cite{Podlubny1999}
\begin{equation}\label{eq.1.2}
\displaystyle
\mathcal{L}^\gamma u(x)=\partial^\gamma_{x,\Omega} u(x):=
-\frac{1}{2\cos\left(\frac{\pi}{2}\gamma\right)\Gamma\left(2-\gamma\right)}
\frac{\mathrm{d}^2}{\mathrm{d}x^2}
\int_{\Omega}|x-s|^{1-\gamma}u\left(s\right)\mathrm{d}s,\;x\in\mathds{R}.
\end{equation}
In the two-dimensional and three-dimensional cases, the Riesz derivative
 can also be defined similarly, that is,
\begin{equation*}
\displaystyle
\mathcal{L}^\gamma u(x,y)=
\left(\partial^\gamma_{x,\Omega}+\partial^\gamma_{y,\Omega}\right)u(x,y),
\;(x,y)\in\mathds{R}^2,
\end{equation*}
and
\begin{equation*}
\displaystyle
\mathcal{L}^\gamma u(x,y,z)=
\left(\partial^\gamma_{x,\Omega}+\partial^\gamma_{y,\Omega}+
\partial^\gamma_{z,\Omega}
\right)u(x,y,z),\;(x,y,z)\in\mathds{R}^3.
\end{equation*}

Similar to the integer-order Allen-Cahn equation,
the space fractional equation (1.1), as a phase-field model,
can also be viewed as an $L^2$ gradient flow with respect
to the following fractional analog of the
Ginzburg-Landau free energy functional

\begin{equation}\label{eq.1.3}
\displaystyle
E\left(u\right)=\int_{\Omega}\left(F\left(u\right)-
\frac{1}{2}\varepsilon^2u\mathcal{L}_{\gamma}u\right)
\mathrm{d}\vecx.
\end{equation}

It is well-known that the exact solution of the Allen-Cahn equation (\ref{eq.1.1})
has two important properties \cite{Hou2017}, i.e.,
preserves the maximum principle
and satisfy the energy dissipation property with time. More precisely, the first one means that
if the initial value is constrained by the constant 1 in the $L^\infty$ norm,
then the entire solution will be constrained by 1, i.e.,
$$\left\|u_0(\vecx)\right\|_{L^\infty}\leq1
\Rightarrow\left\|u(\vecx,t)\right\|_{L^\infty}\leq1
\;\mathrm{for}\;\mathrm{all}\;(\vecx,t)\in\overline{\Omega}\times(0,T],$$
where the symbol $\left\|\cdot\right\|_{L^\infty}$ is the $L^\infty$ norm and
defined by
$
\left\|u\right\|_{L^\infty}=\max_{\vecx\in\overline{\Omega}}
\left|u(\vecx)\right|.$
The second one implies that the corresponding
energy functional (\ref{eq.1.3}) decreases over time, that is
\begin{equation*}
\displaystyle
\frac{\mathrm{d}}{\mathrm{d}t}E\left(u\right)=
\left(\frac{\delta E(u)}{\delta u},\phi(u)\right)
=-\left\|\partial_t\left(\vecx,t\right)\right\|^2
\leq0,\;\mathrm{for}\;\mathrm{any}\;t>0,
\end{equation*}
which is often called the nonlinear energy stability,
where $\phi(u)=\epsilon^2\mathscr{L}^\gamma u+f(u)$,
the symbols $(\cdot, \cdot)$ and $\left\|\cdot\right\|$ are the standard inner
product and $L^2$ norm, which are defined by
\begin{equation*}
\displaystyle
\left(u,v\right)=\int_{\Omega}uv\mathrm{d}\vecx,\;
\;\;
\left\|u\right\|^2=\int_{\Omega}\left|u\right|^2\mathrm{d}\vecx,\;\;
\mathrm{any}\;u,v\in L^2(\Omega).
\end{equation*}

Spatial nonlocal Allen-Cahn equation (\ref{eq.1.1}) is
a partial differential equation, which extends the classical
 Allen-Cahn equation by introducing fractional derivative. The classic Allen-Cahn equation is one of the basic
  models to study phase transition and interface dynamics,
  especially material science and fluid dynamics. Because
  the nonlocal model introduces the fractional derivative
  term considering the long-range interaction, it is suitable
   for describing a wider range of phenomena.
For example, nonlocal Allen-Cahn equation (\ref{eq.1.1}) can be
used to simulate the dynamic behavior of cell populations or
tissues, in which cell interactions over long distances can be
 represented by fractional derivatives. For example, in the
 tumor growth model, the fractional derivative can explain the
  influence of nutrient concentration or signal molecules on
  diffusion in a large range. In addition, this model can also
   be used to study the pattern formation process in biological
    tissues, such as the arrangement and combination of cells
    in developing embryos, in which long-range interaction
    plays a crucial role, that is, fractional derivatives
     play an important role.
Unfortunately, due to the influence of nonlocal operators,
the exact solution of model (\ref{eq.1.1}) is generally difficult to
obtain, which has caused great obstacles to the further
study of the dynamic behavior of the equation, so it
is natural to seek an efficient numerical method as one of the best ways.

From equation (\ref{eq.1.1}), it is not difficult to find that
in order to derive its efficient numerical algorithm,
the first task is to construct the higher-order
numerical differential formula of Riesz derivative.
Generally speaking, there are two ways to construct
the numerical differential formula of Riesz derivative,
indirect method and direct method. In more detail,
the former is to construct the numerical differential
 formulas of the left and right Riemann-Liouville derivatives
 first, and then sum them by weight. The latter is constructed
 from the definition itself.
Up to now, some scholars have constructed a series of
 efficient numerical differential formulas in these two ways.
In the direct method, the Gr\"{u}nwald–Letnikov formulas was used
to approximate the Riemann-liouville derivatives in the early stage.
Unfortunately, when it is applied to time-dependent spatial
fractional differential equations, the numerical algorithm
obtained is unstable. In order to make up for this
shortcoming, Meerschaert and Tadjeran \cite{Meerschaert2004} proposed the
so-called shifted Gr\"{u}nwald–Letnikov formula with
convergence order of 1.
Subsequently, Deng's group put forward second-order
 and third-order numerical differential formulas
 by weighting and shifting Gr\"{u}nwald–Letnikov formulas
 \cite{Tian2015,Zhou2013}.
In order to improve stability and accuracy, Ding \cite{Ding2017} also
constructed a second-order numerical differential formula by
 constructing a new generation function.
Looking back at the direct method, the most famous one is the
second-order formula \cite{Celik2012}based on the fractional central
difference operator proposed by Ortigueira \cite{Ortigueira2006}.
In recent years, higher-order formulas with this second-order
formula as the core have also been obtained
\cite{Ding2015,Zhao2014,Ding20172}.
But as far as we know, numerical differential formulas
with convergence order of 6 are still scarce.

Based on these numerical differential formulas mentioned
 above and their improvements and corrections, a
 large number of numerical algorithms have emerged
  for phase field models (\ref{eq.1.1}). Such as in \cite{Hou2017},
Hou et al. used the Crank-Nicolson method in time and second-order
numerical differential formula by developed Tian et al.
\cite{Tian2015} in space,
and constructed a difference scheme with order
$\mathcal{O}\left(\tau^2+h^2\right)$, which show that
the numerical solution satisfied discrete maximum principle under
$0<\tau\leq\min \left\{\frac{1}{2},\frac{h^\gamma}{2d\epsilon^2}\right\}$.
Using the fourth-order numerical differential formula proposed by
Ding and Li in \cite{Ding2015},
He et al. \cite{He2020} constructed
a spatial fourth-order operator splitting
scheme for equation (\ref{eq.1.1}),
and the discrete maximum principle was discussed
under time step $\tau$ satisfies
certain conditions.
However, it is a pity that they do not give detailed theoretical proof for the
error estimation and energy stability of the numerical scheme.
In addition, from the numerical experiments given by them,
it is not difficult to find that the energy is not always reduced, which reflects
 that the numerical scheme established does
 not always satisfy the energy dissipation property.
In \cite{Zhang2021}, Zhang
et al. presented an implicit-explicit Runge-Kutta schemes with order
$\mathcal{O}\left(\tau^p+h^2\right)$, where $p\geq1$.
And the discrete maximum principle,
energy stability and error estimation were studied
under the assumption of time step $0<\tau\leq\frac{h^\gamma}{dg_0\epsilon^2}\tilde{r}$,
 where $g_0=\frac{\Gamma(\gamma+1)}{\Gamma^2(\gamma/2+1)}$ and $\tilde{r}$ is a positive
  constant satisfying certain conditions. However, it is a little regrettable that the theoretical
   analysis of energy dissipation is not given in a strict sense, they only give
   the boundedness of discrete energy, which shows that the energy
   dissipation of the numerical solution is inconsistent with that of the exact solution.
More recently, with the help of
the symplectic Runge-Kutta method and the quadratic
auxiliary variable technique \cite{Gong2022} in time, and
a fourth-order numerical differential formula developed by Hao et al. in \cite{Hao2021}
 is applied to the space Riesz derivative.
Xu and Fu \cite{Xu2024} proposed a fully discretized implicit scheme with order
$\mathcal{O}\left(\tau^2+h^4\right)$,  they analyzed the discrete energy stability and
maximum principle under the condition
$0\leq\tau\leq\min\left\{\frac{1}{2+\delta},\frac{1}{3},\frac{h^\gamma}{2d\epsilon^2}\right\}$,
where $0\leq\delta\leq4$. For more studies on equation (\ref{eq.1.1}),
it can be found in references \cite{Bu2019,Chen12021,Chen22021}.

From above analysis, we can see that although
there are some studies on the space fractional Allen-Cahn equation (\ref{eq.1.1}),
it is still limited.
Moreover, as far as we know, there is no numerical scheme whose spatial
convergence order is higher than 4, and the condition of preserving
the discrete maximum principle
is independent of $\epsilon$ and $h$.
Therefore, the main motivation of our work in this paper is to
 make up for these shortcomings, so as to obtain an efficient
 numerical scheme for solving equation (\ref{eq.1.1}).
Compared with existing researches, our contributions are listed as below.
\begin{itemize}
\item By constructing a new generating function, a sixth-order numerical differential
 formula for approximating the space Riesz derivative is established.

\item The implicit integral factor method coupled with Pad\'{e} approximation strategy
is proposed to deal with the time derivative.

\item The condition of preserving the maximum principle of the numerical scheme constructed
in this paper is $0<\tau\leq1$, which is independent of the parameters
$\epsilon$ and $h$, and
 significantly improves the condition of preserving the maximum principle
  of the numerical schemes in the existing literatures.

\item The convergence order of our difference scheme can reach
$\mathcal{O}\left(\tau^2+h^6\right)$
 in the sense of discrete maximum norm, which is obviously higher
 than the existing schemes.
\end{itemize}

The rest of this paper is organized as follows. In
Section 2, we construct a new sixth-order numerical differential
 formula that approximates
the Riesz derivative, and apply it to equation (\ref{eq.1.1}),
thus obtaining a fully discrete finite
difference scheme.
In Section 3, the properties of the proposed difference
scheme are studied in detail,
 including discrete maximum principle, discrete energy
 stability and convergence.
 Numerical experiments are presented
to confirm the theoretical analysis and demonstrate the
efficiency of the provided
schemes in Section 4. The last section is concerned with the conclusion.

\section{Numerical scheme}
In this section, we first present a sixth-order numerical differential formula
for the space Riesz derivative, and obtain the corresponding semi-discrete numerical scheme.
Then a fully discrete scheme is proposed by applying the
implicit integration factor method couple Pad\'{e} approximation strategy in time.
\subsection{Space discretization}Here, we apply a uniform
 partition of $\Omega$ into square elements of length $h$,
and let $M$
be the total number of grid nodes.
In order to effectively approximate the Riesz derivative, we can construct a higher-order
numerical differential formula by using the strategy
 of seeking suitable new generating functions. In the paper, we design a new
 generating function $G(w)$ with the following form
\begin{equation}\label{eq.2.1}
\begin{aligned}\displaystyle
{G}(w)=\left(2-w-w^{-1}\right)^{\frac{\gamma}{2}}\left[1+\frac{\gamma}{24}
\left(2-w-w^{-1}\right)+\frac{\gamma\left(5\gamma+22\right)}{5760}
 \left(2-w-w^{-1}\right)^2\right],
\end{aligned}
\end{equation}
and it is analytic, that is
\begin{equation}\label{eq.2.2}
\begin{aligned}\displaystyle
\displaystyle{G}(w)=
\sum\limits_{m=-\infty}^{+\infty}
g_{m}^{(\gamma)}w^m,\;\;|w|\leq1.
\end{aligned}
\end{equation}
Using a similar technique as in \cite{Ding2022,Ding2023}, we know that the exact
expression of the coefficient $g_{m}^{(\gamma)}$ is
\begin{equation}\label{eq.2.3}
\begin{aligned}\displaystyle
g_{m}^{(\gamma)}=
&\displaystyle\frac{1}{2\pi}\int^{\pi}_{-\pi}
{G}\left(e^{\mathrm{i}s}\right)
e^{-\mathrm{i}ms}\mathrm{d}s\vspace{0.4cm}\\
=&\displaystyle
\frac{(-1)^m \Gamma\left(\gamma+1\right)}
{\Gamma\left(\frac{\gamma}{2}-m+1\right)\Gamma\left(\frac{\gamma}{2}+m+1\right)}
\left[1+\frac{\gamma(\gamma+1)(\gamma+2)}{6\left(\gamma-2m+2\right)
\left(\gamma+2m+2\right)}\right.\vspace{0.2cm}\\
&\displaystyle \left.+\frac{\gamma(\gamma+1)(\gamma+2)(\gamma+3)(\gamma+4)(5\gamma+22)}
{360\left(\gamma-2m+4\right)\left(\gamma-2m+2\right)
\left(\gamma+2m+4\right)\left(\gamma+2m+2\right)}
\right],\;m=0,\pm1,\cdots.
\end{aligned}
\end{equation}

Next, based on the above generating function (\ref{eq.2.1}),
a new high-order
numerical differential formula is constructed
for approximating the Riesz derivative (\ref{eq.1.2}), and the main result is
described in the following theorem.

\begin{theorem}\label{Th.2.1}
Suppose
\begin{equation*}
\begin{aligned}\displaystyle
u\in\mathscr{C}^{{6+\gamma}}\left(\mathds{{R}}\right)=\left\{u\,|\,u\in L^1(\mathds{{R}}),
\;\mathrm{and}\;\int_{\mathds{{R}}}\left(1+|s|\right)^{{6+\gamma}}\,|\,\hat{u}(s)
\,|\,\mathrm{d}s<\infty\right\},\;\;1<\gamma\leq2,
\end{aligned}
\end{equation*}
where $\hat{u}(s)$ is the Fourier transform of $u(z)$ for all $s\in \mathds{{R}}$.
Then there holds that
\begin{equation*}
\begin{aligned}\displaystyle
\partial_{x,\mathds{R}}^\gamma u(x)= {\delta}_{h,\mathds{R}}^{\gamma}u(x)
+\mathcal{O}\left(h^6\right)
\end{aligned}
\end{equation*}
uniformly as $h\rightarrow0$, where the fractional difference operator
${\delta}_{h,\mathds{R}}^{\gamma}$ is defined as
\begin{equation*}
\begin{aligned}
\displaystyle {\delta}_{h,\mathds{R}}^{\gamma}u(x)
=-{h^{-\gamma}}\sum\limits_{m=-\infty}^{+\infty}
g_{m}^{(\gamma)}u\left(x-mh\right).
\end{aligned}
\end{equation*}
\end{theorem}

\begin{proof}
Denote
\begin{equation}\label{eq.2.4}
\begin{aligned}\displaystyle
\Phi(\varpi,h)=\mathcal{F}\left\{\partial_{x,\mathds{R}}^\gamma u(x)\right\}-
\mathcal{F}\left\{{\delta}_{h,\mathds{R}}^{\gamma}u(x)\right\}.
\end{aligned}
\end{equation}

Note that the Fourier transforms of the operators $\partial_{x,\mathds{R}}^\gamma$
and ${\delta}_{h,\mathds{R}}^{\gamma}$ are
\begin{equation*}
\begin{aligned}\displaystyle
\mathcal{F}\left\{\partial_{x,\mathds{R}}^\gamma\right\}
=-|\varpi|^\gamma\hat{u}(\varpi),
\end{aligned}
\end{equation*}
and
\begin{equation*}
\begin{aligned}\displaystyle
\mathcal{F}\left\{{\delta}_{h,\mathds{R}}^{\gamma}u(x)\right\}
=&-{h^{-\gamma}}\sum\limits_{m=-\infty}^{+\infty}
g_{m}^{(\gamma)}\mathcal{F}\left\{u\left(x-mh\right),\varpi\right\}=
-{h^{-\gamma}}\sum\limits_{m=-\infty}^{+\infty}
g_{m}^{(\gamma)}e^{\mathrm{i}m\varpi h}\hat{u}(\varpi)\\
 =&-|\varpi|^\gamma\left[1+\frac{\gamma}{6}\sin^2\left(\frac{1}{2}\varpi h\right)
 +\frac{\gamma(5\gamma+22)}{360}\sin^4\left(\frac{1}{2}\varpi h\right) \right]
 \left|\frac{\sin\left(\frac{1}{2}\varpi h\right)}{\frac{1}{2}\varpi h}\right|^\gamma\hat{u}(\varpi),
\end{aligned}
\end{equation*}
respectively, which further leads to
\begin{equation*}
\begin{aligned}\displaystyle
\Phi(\varpi,h)=|\varpi|^\gamma
\left\{-1+\left[1+\frac{\gamma}{6}\sin^2\left(\frac{1}{2}\varpi h\right)
 +\frac{\gamma(5\gamma+22)}{360}\sin^4\left(\frac{1}{2}\varpi h\right) \right]
 \left|\frac{\sin\left(\frac{1}{2}\varpi h\right)}{\frac{1}{2}\varpi h}\right|^\gamma\right\}\hat{u}(\varpi).
\end{aligned}
\end{equation*}

From the Taylor expansions, there are
\begin{equation*}
\begin{aligned}\displaystyle
1+\frac{\gamma}{6}\sin^2\left(\frac{1}{2}\varpi h\right)
 +\frac{\gamma(5\gamma+22)}{360}\sin^4\left(\frac{1}{2}\varpi h\right)
 =1+\frac{\gamma}{24}\left(\varpi h\right)^2+\frac{\gamma(5\gamma+2)}{5760}
 \left(\varpi h\right)^4+\mathcal{O}\left(\varpi h\right)^6,
\end{aligned}
\end{equation*}
and
\begin{equation*}
\begin{aligned}\displaystyle
 \left|\frac{\sin\left(\frac{1}{2}\varpi h\right)}{\frac{1}{2}\varpi h}\right|^\gamma
 =1-\frac{\gamma}{24}\left(\varpi h\right)^2+\frac{\gamma(5\gamma-2)}{5760}
 \left(\varpi h\right)^4+\mathcal{O}\left(\varpi h\right)^6,
\end{aligned}
\end{equation*}
we finally know that
\begin{equation*}
\begin{aligned}\displaystyle
\displaystyle \Phi(\varpi,h)=|\varpi|^\gamma\left[-\frac{\gamma^2}{34560}
\left(\varpi h\right)^6
+\frac{\gamma^2\left(25\gamma^2-4\right)}{33177600}\left(\varpi h\right)^8
+\mathcal{O}\left(\varpi h\right)^{10}
\right]\hat{u}(\varpi),
\end{aligned}
\end{equation*}
and there exist a positive constant $C_1$, such that
\begin{equation}\label{eq.2.5}
\begin{aligned}\displaystyle
\displaystyle \left|\Phi\left(\varpi,h\right)\right|
\leq C_1\left(\left|\varpi\right|^{\gamma+6}\left|\hat{u}\left(\varpi\right)\right|\right)h^{6}.
\end{aligned}
\end{equation}

Using the inverse Fourier transform to equation (\ref{eq.2.4}), and
combining with inequality (\ref{eq.2.5})
 and condition $u(x)\in \mathscr{C}^{6+\gamma}(\mathds{{R}})$,
 we can finally achieve the following result
\begin{equation*}
\begin{aligned}
\displaystyle \left|
\partial_{x,\mathds{R}}^\gamma u(x)-
{\delta}_{h,\mathds{R}}^{\gamma}u(x)\right|
=&\displaystyle\left|\frac{1}{2\pi}\int_{{\mathds{R}}}\Phi(\varpi,h)
e^{\mathrm{i}\varpi x}\mathrm{d}\varpi\right|\leq
\frac{1}{2\pi}\int_{{\mathds{R}}}\left|\Phi(\varpi,h)\right|
\mathrm{d}\varpi\vspace{0.1 cm}\\
 \leq
&\displaystyle\frac{C_1}{2\pi}\left(\int_{{\mathds{R}}}|\varpi|^{\gamma+6}|\hat{u}(\varpi)|
\mathrm{d}\varpi\right)h^{6}
\leq \frac{C_1}{2\pi}\left(\int_{{\mathds{R}}}\left(1+|\varpi|^{\gamma+6}\right)|\hat{u}(\varpi)|
\mathrm{d}\varpi\right)h^{6}\\
=&\displaystyle C_2h^6.
\end{aligned}
\end{equation*}
Obviously, this is completely consistent with the result to be proved, which completes the proof.
\end{proof}

Now, we consider the case for $x\in \Omega=[a,b]\subset\mathds{{R}}$.
If $u^{*}(x)$ is defined by
\begin{eqnarray*}
u^{*}(x)=\left\{
\begin{array}{lll}
u(x),\;\;x\in\Omega,\\
0,\;\;x\notin\Omega,
\end{array}
\right.
\end{eqnarray*}
such that $u^{*}(x)\in
\mathscr{C}^{6+\gamma}(\mathds{{R}})$ , then it follows from the Theorem \ref{Th.2.1} that
$
\partial_{x,\mathds{R}}^\gamma u^{*}(x)= {\delta}_{h,\mathds{R}}^{\gamma}u^{*}(x)
+\mathcal{O}\left(h^6\right).
$
Because of $u^{*}(x)=0$ for $x\notin\Omega$, then we further obtain
\begin{equation}\label{eq.2.6}
\begin{aligned}\displaystyle
\partial_{x,\Omega}^\gamma u(x)= {\delta}_{h,\Omega}^{\gamma}u(x)
+\mathcal{O}\left(h^6\right),
\end{aligned}
\end{equation}
where the operator
${\delta}_{h,\Omega}^{\gamma}$ is given by
\begin{equation}\label{eq.2.7}
\begin{aligned}\displaystyle
{\delta}_{h,\Omega}^{\gamma}u(x)=-{h^{-\gamma}}
\sum\limits_{m=(x-b)/h}^{(x-a)/h}g_{m}^{(\gamma)}u\left(x-mh\right).
\end{aligned}
\end{equation}

 Let $M$ be a positive integer, and
 $x_j=y_j=z_j=a+jh,j=0,1,\cdots,M$, with $h=(b-a)/M$.
 Therefore, under homogeneous Dirichlet boundary conditions, it follows from
 equations (\ref{eq.2.6}) and (\ref{eq.2.7}), we
 can see that the discretization
matrix of  the operator $\mathcal{L}^\gamma$ with mesh size $h$
 in one-dimensional is given by
 \begin{equation*}
\begin{aligned}\displaystyle
\mathbf{D}_\gamma^{(1)}=-\frac{1}{h^\gamma}\mathbf{K}_\gamma,
\end{aligned}
\end{equation*}
 with
$$\mathbf{K}_\gamma=
\left(
  \begin{array}{ccccc}
   g_{0}^{(\gamma)} & g_{1}^{(\gamma)}& g_{2}^{(\gamma)}
    & \cdots & g_{M-2}^{(\alpha)} \vspace{0.2 cm}\\
   g_{-1}^{(\gamma)} & g_{0}^{(\gamma)} & g_{1}^{(\gamma)}
    & g_{2}^{(\gamma)}& \cdots \vspace{0.2 cm}\\
    \vdots &\vdots & \ddots & \ddots & \ddots \vspace{0.2 cm}\\
  g_{3-M}^{(\gamma)} & \ldots & g_{-1}^{(\alpha)}&g_{0}^{(\gamma)}&
  g_{1}^{(\gamma)} \vspace{0.2 cm}\\
g_{2-M}^{(\gamma)}& g_{3-M}^{(\gamma)}& \ldots & g_{-1}^{(\gamma)}&
  g_{0}^{(\gamma)} \vspace{0.2 cm}\\
  \end{array}
\right)_{(M-1)\times(M-1)}.
$$

Similarly, it is easy to derive
the corresponding
discretization matrices of the operator $\mathcal{L}^\gamma$ in two-dimensions
and in three-dimensions
are
$$\mathbf{D}_\gamma^{(2)}=-\frac{1}{h^\gamma}\left(\mathbf{I}\otimes\mathbf{K}_\gamma
+\mathbf{K}_\gamma\otimes \mathbf{I}\right),$$
and$$
\mathbf{D}_\gamma^{(3)}=-\frac{1}{h^\gamma}\left(\mathbf{I}\otimes \mathbf{I}
\otimes\mathbf{K}_\gamma
+\mathbf{I}\otimes\mathbf{K}_\gamma\otimes \mathbf{I}
+\mathbf{K}_\gamma\otimes \mathbf{I}\otimes \mathbf{I}\right),$$
respectively, where $\mathbf{I}$ is an identity
 matrix of size $M-1$ and the symbol $\otimes$ represents
the Kronecker product between the two matrices.

Once again, we study the properties of the coefficients
$g_{m}^{(\gamma)}(m=0,\pm1,\pm2,\cdots)$. The main results are
 described in the following theorem.

\begin{lemma}\label{Lem.2.2} For the coefficients
$g_{m}^{(\gamma)}(m=0,\pm1,\pm2,\cdots)$, the following properties hold:\vspace{0.2cm}\\
$\mathrm{(1)}$~~$\displaystyle g_{0}^{(\gamma)}\geq0$,\;
$ g_{\pm1}^{(\gamma)}\leq0\;$;\vspace{0.2 cm}\\
$\mathrm{(2)}$~~$g_{\pm2}^{(\gamma)}\leq0\; \mathrm{if}\;\gamma\in(1,\gamma^{*}],
\;\mathrm{while} \;g_{\pm2}^{(\gamma)}\geq0\; \mathrm{if}\;\gamma\in[\gamma^{*},2).
\;\mathrm{Here}$,\;
$\gamma^{*}=1.4746120;$\vspace{0.2 cm}\\
$\mathrm{(3)}$~~$ g_{m}^{(\gamma)}\leq0,\;m=\pm3,\pm4,\cdots$;\vspace{0.2 cm}\\
$\mathrm{(4)}$~~$\displaystyle\sum\limits_{m=-\infty}^{+\infty}g_{m}^{(\gamma)}=0.$
\end{lemma}

\begin{proof}
(1)~From equation (\ref{eq.2.3}), we know that
\begin{equation*}
\begin{aligned}\displaystyle
g_{0}^{(\gamma)}=
\frac{\left(\gamma+2\right)\left(\gamma+4\right)\left(5\gamma^4
+102\gamma^3+763\gamma^2
+2466\gamma+2880\right)
\Gamma\left(\gamma+1\right)}{5760\Gamma^2\left(\gamma/2+3\right)}>0,
\end{aligned}
\end{equation*}
and that
\begin{equation*}
\begin{aligned}\displaystyle
g_{\pm1}^{(\gamma)}=-
\frac{\left(5\gamma^4+122\gamma^3+1171\gamma^2
+5278\gamma+9624\right)
\Gamma\left(\gamma+1\right)}{90\gamma\left(\gamma+2\right)
\left(\gamma+4\right)\left(\gamma+6\right)
\Gamma^2\left(\gamma/2\right)}<0.
\end{aligned}
\end{equation*}

$\mathrm{(2)}$~In equation (\ref{eq.2.3}), we take $m=\pm2$, and denote
\begin{equation}\label{eq.2.8}
\begin{aligned}\displaystyle
g_{\pm2}^{(\gamma)}=p_1\left(\gamma\right)q_1\left(\gamma\right),
\end{aligned}
\end{equation}
where
\begin{equation}\label{eq.2.9}
\begin{aligned}\displaystyle
p_1\left(\gamma\right)=5\gamma^5+132\gamma^4+1415\gamma^3
+6900\gamma^2+9380\gamma-34032,
\end{aligned}
\end{equation}
and
\begin{equation}\label{eq.2.10}
\begin{aligned}\displaystyle
q_1\left(\gamma\right)=
\frac{\Gamma\left(\gamma+1\right)}{90\gamma\left(\gamma+2\right)
\left(\gamma+4\right)\left(\gamma+6\right)\left(\gamma+8\right)
\Gamma^2\left(\gamma/2\right)}>0.
\end{aligned}
\end{equation}

By Newton iteration method, we get a unique solution of the equation $p_1\left(\gamma\right)=0$
for $\gamma\in(1,2]$ is
$\gamma^{*}=1.474612$. Combining (\ref{eq.2.8}), (\ref{eq.2.9}) and (\ref{eq.2.10}),
 we can see that the conclusion is valid.

$\mathrm{(3)}$~Denote
$
g_{m}^{(\gamma)}=p_2(\gamma)\left[q_2(\gamma)+q_3(\gamma)\right],
$
where
\begin{equation*}
\begin{aligned}\displaystyle
p_2(\gamma)=
\frac{(-1)^m \Gamma\left(\gamma+1\right)}
{\Gamma\left({\gamma}/2-m+1\right)\Gamma\left({\gamma}/2+m+1\right)},\;\;
q_2(\gamma)=1+\frac{\gamma\left(\gamma+1\right)\left(\gamma+2\right)}
{6\left(\gamma-2m+2\right)
\left(\gamma+2m+2\right)},
\end{aligned}
\end{equation*}
and
\begin{equation*}
\begin{aligned}\displaystyle
q_3(\gamma)=
\frac{\gamma(\gamma+1)(\gamma+2)(\gamma+3)(\gamma+4)(5\gamma+22)}
{360\left(\gamma-2m+4\right)\left(\gamma-2m+2\right)
\left(\gamma+2m+4\right)\left(\gamma+2m+2\right)}.
\end{aligned}
\end{equation*}

From \cite{Celik2012}, we know that $p_2(\gamma)\leq0$ for all $|m|\geq1$ and $\gamma\in(1,2]$.
Meanwhile, for $\gamma\in(1,2]$ and all $|m|\geq3$,
there also holds that
\begin{equation*}
\begin{aligned}\displaystyle
q_2(\gamma)=&\displaystyle1+\frac{\gamma(\gamma+1)(\gamma+2)}
{6\left(\gamma-2m+2\right)
\left(\gamma+2m+2\right)}\geq
\left[1+\frac{\gamma(\gamma+1)(\gamma+2)}{6\left(\gamma-2m+2\right)
\left(\gamma+2m+2\right)}\right]\Bigg|_{m=\pm3}\\\displaystyle
=&\displaystyle\frac{144-\gamma(\gamma^2+9\gamma+14)}{6\left(\gamma+6\right)
\left(4-\gamma\right)}>0,
\end{aligned}
\end{equation*}
and
\begin{equation*}
\begin{aligned}\displaystyle
q_3(\gamma)>
\frac{\gamma(\gamma+1)(\gamma+2)(\gamma+3)(\gamma+4)(5\gamma+22)}
{360}
\left(\frac{1}{2m}\right)^4>0.
\end{aligned}
\end{equation*}
This shows that
$g_{m}^{(\gamma)}=p_2(\gamma)\left[q_2(\gamma)+q_3(\gamma)\right]\leq0$
for all $|m|\geq3$ and $\gamma\in(1,2]$.

$\mathrm{(4)}$~By taking $w=1$ in equation (\ref{eq.2.2}) and combining it
with equation (\ref{eq.2.1}), the result to be
proved can be obtained immediately, which completes all the proofs.
\end{proof}

After studying the high-order numerical differential formula of Riesz derivative and its
 corresponding coefficients, we immediately turn to consider the equation (\ref{eq.1.1}).
Denote $$\mathbf{u}(\vecx,t)=\left[u(\vecx_1,t),u(\vecx_2,t),
\cdots,u(\vecx_{M-1},t)\right]^{\mathsf{T}}$$ as
the solution vector of (\ref{eq.1.1}) on the spatial mesh, using the numerical differential
 formula (\ref{eq.2.6})
to the equation (\ref{eq.1.1}), then we
obtain the corresponding semi-discrete scheme as follows
\begin{eqnarray}\label{eq.2.11}
\begin{array}{lll}\displaystyle
\frac{\mathrm{d}\mathbf{u}(\vecx,t)}{\mathrm{d}t}+\mathbf{A}^{(d)}\mathbf{u}(\vecx,t)
=\mathbf{f}\left(\mathbf{u}\right)+R_h,\;d=1,2,3,
\end{array}
\end{eqnarray}
where $\mathbf{A}^{(d)}=-\epsilon^2\mathbf{D}_\gamma^{(d)}$, and
$R_h$ is the local truncation error in spatial direction.
Meanwhile,
it follows from the numerical differential
 formula (\ref{eq.2.6}), there holds that
\begin{eqnarray}\label{eq.2.12}
\begin{array}{lll}\displaystyle
\left\|R_h\right\|_\infty\leq C_3h^6
\end{array}
\end{eqnarray}
as the exact solution $\mathbf{u}(\vecx,\cdot)\in\mathscr{C}^{6+\gamma}
\left(\overline{\Omega}\right)$,
where $C_3$ is a positive constant depending on the exact
 solution $\mathbf{u}(\vecx,t)$
but not depending on $h$.

\subsection{Time discretization}
For the time discretization, let us consider a uniform partition of the time interval
$[0, T ] : 0 = t_0 < t_1 <\cdots < t_N = T $, with the step size $\tau= T/N$.
Using implicit integration factor method \cite{Nie2008}, the exact solution of the
 equation (\ref{eq.2.11})
is given by the following equation
\begin{eqnarray*}
\begin{array}{lll}\displaystyle
\mathbf{u}\left(\vecx,t\right)=e^{-t\mathbf{A}^{(d)}}u_0+\int_{0}^{t}
e^{-(t-s)\mathbf{A}^{(d)}}\mathbf{f}\left(\mathbf{u}\left(\vecx,s\right)\right)\mathrm{d}s
+\tau R_h,
\end{array}
\end{eqnarray*}
which satisfies the recurrence formula
\begin{eqnarray}\label{eq.2.13}
\begin{array}{lll}\displaystyle
\mathbf{u}\left(\vecx,t_{k+1}\right)=e^{-\tau \mathbf{A}^{(d)}}\mathbf{u}
\left(\vecx,t_{k}\right)+\int_{0}^{\tau}
e^{-(\tau-s)\mathbf{A}^{(d)}}\mathbf{f}\left(\mathbf{u}
\left(\vecx,t_{k}+s\right)\right)\mathrm{d}s+
\tau R_h.
\end{array}
\end{eqnarray}

By using the [1,1] pad\'{e} approximation \cite{Moler2003} to
deal with the exponential matrices and
the trapezoidal method to compute
the resulting integral in equation (\ref{eq.2.13}), we obtain the
 following equation
\begin{equation}\label{eq.2.14}
\begin{aligned}\displaystyle
\mathbf{u}\left(\vecx,t_{k+1}\right)=&\left(
\mathbf{I}+\frac{\tau}{2}\mathbf{A}^{(d)}
\right)^{-1}\left(
\mathbf{I}-\frac{\tau}{2}\mathbf{A}^{(d)}
\right)
\mathbf{u}\left(\vecx,t_{k}\right)
\\&+\int_{0}^{\tau}
\left(
\mathbf{I}+\frac{\tau-s}{2}\mathbf{A}^{(d)}
\right)^{-1}\left(
\mathbf{I}-\frac{\tau-s}{2}\mathbf{A}^{(d)}
\right)
\mathbf{f}\left(\mathbf{u}\left(\vecx,t_{k}+s\right)\right)\mathrm{d}s
\\&+\tau \left(R_\tau+ R_h\right)\\
=&\left(
\mathbf{I}+\frac{\tau}{2}\mathbf{A}^{(d)}
\right)^{-1}\left(
\mathbf{I}-\frac{\tau}{2}\mathbf{A}^{(d)}
\right)\mathbf{u}\left(\vecx,t_{k}\right)
\\&+\frac{\tau}{2}\mathbf{f}\left(\mathbf{u}\left(\vecx,t_{k+1}\right)\right)+
\frac{\tau}{2}\left(
\mathbf{I}+\frac{\tau}{2}\mathbf{A}^{(d)}
\right)^{-1}\left(
\mathbf{I}-\frac{\tau}{2}\mathbf{A}^{(d)}
\right)
\mathbf{f}\left(\mathbf{u}\left(\vecx,t_{k}\right)\right)\\&
+\tau \left(R_\tau+ R_h\right),
\end{aligned}
\end{equation}
where $R_\tau$ is the local truncation error in temporal direction and satisfies
\begin{eqnarray}\label{eq.2.15}
\begin{array}{lll}\displaystyle
\left\|R_\tau\right\|_\infty\leq C_4\tau^2
\end{array}
\end{eqnarray}
for $\mathbf{u}(\cdot,t)\in{C}^{2}\left[0,T\right]$,
where $C_4$ is also a positive constant that depends on the exact
solution $\mathbf{u}(\vecx,t)$ rather than $\tau$.

Omitting the truncation error $\tau \left(R_\tau+ R_h\right)$
 in the above equation, and replacing exact solution $\mathbf{u}\left(\vecx,t_{k}\right)$ by its
numerical solution $\mathbf{U}^{k}$, then we can obtain
a high-order difference scheme for equation (\ref{eq.1.1}) as follows
\begin{equation}\label{eq.2.16}
\begin{aligned}\displaystyle
\mathbf{U}^{k+1}=&\left(
\mathbf{I}+\frac{\tau}{2}\mathbf{A}^{(d)}
\right)^{-1}\left(
\mathbf{I}-\frac{\tau}{2}\mathbf{A}^{(d)}
\right)\mathbf{U}^{k}+
\frac{\tau}{2}\mathbf{f}
\left(\mathbf{U}^{k+1}\right)\\&+
\frac{\tau}{2}\left(
\mathbf{I}+\frac{\tau}{2}\mathbf{A}^{(d)}
\right)^{-1}\left(
\mathbf{I}-\frac{\tau}{2}\mathbf{A}^{(d)}
\right)\mathbf{f}\left(\mathbf{U}^{k}\right)
,\;\;k=0,1,\cdots,N-1.
\end{aligned}
\end{equation}

\section{Theoretical analysis of the numerical scheme}
In this section, for simplicity, we only consider the one dimensional
space fractional Allen-Cahn equation, i.e., $d=1$.
From equation (\ref{eq.2.16}), it is not difficult to obtain the corresponding numerical scheme
for case $d=1$ as follows
\begin{equation}\label{eq.3.1}
\begin{aligned}\displaystyle
\mathbf{U}^{k+1}=&\left(
\mathbf{I}+\frac{\tau}{2}\mathbf{A}^{(1)}
\right)^{-1}\left(
\mathbf{I}-\frac{\tau}{2}\mathbf{A}^{(1)}
\right)\mathbf{U}^{k}+
\frac{\tau}{2}\mathbf{f}
\left(\mathbf{U}^{k+1}\right)\\&+
\frac{\tau}{2}\left(
\mathbf{I}+\frac{\tau}{2}\mathbf{A}^{(1)}
\right)^{-1}\left(
\mathbf{I}-\frac{\tau}{2}\mathbf{A}^{(1)}
\right)\mathbf{f}\left(\mathbf{U}^{k}\right)
,\;\;k=0,1,\cdots,N-1.
\end{aligned}
\end{equation}

\subsection{Preliminaries}
For the convenience of subsequent analysis, we list some basic symbols and related lemmas in this section.
First, define the following space
\begin{eqnarray*}
\begin{array}{lll}\displaystyle
\mathring{\Psi}_h=\left\{\mathbf{U}:\mathbf{U}=\left\{U_j\right\}\;
\mathrm{is\; grid\; function\; and\;} U_0=U_M=0\right\},
\end{array}
\end{eqnarray*}
equipped with discrete inner product and the associated $L^2$
norm defined as
\begin{eqnarray*}
\begin{array}{lll}\displaystyle
\left\langle\mathbf{U},\mathbf{V}\right\rangle:=h\mathbf{V}^{\mathsf{T}}\mathbf{U}
=h\sum_{j=1}^{M-1}U_jV_j,\;\;
\left\|\mathbf{U}\right\|_2=\sqrt{\left\langle\mathbf{U},\mathbf{U}\right\rangle}\;,\;\mathrm{for}\;
\mathrm{any}\; \mathbf{U},\mathbf{V}\in\mathring{\Psi}_h.
\end{array}
\end{eqnarray*}
Besides, we can also give the following discrete $L^\infty$ norm
\begin{eqnarray*}
\begin{array}{lll}\displaystyle
\left\|\mathbf{U}\right\|_\infty=\max_{1\leq j\leq M-1}\left|U_j\right|.
\end{array}
\end{eqnarray*}

\begin{definition}(See \cite{Chan2007})\label{Def:1}
Let $n\times n$ Toeplitz matrix {$\mathbf{Q}_n$} be in the form:
$$\mathbf{Q}_n=
\left(
  \begin{array}{ccccc}
    q_0 & q_{-1} & \cdots & q_{2-n} & q_{1-n} \vspace{0.2 cm}\\
    q_{1} & q_{0} & q_{-1} & \cdots & q_{2-n}\vspace{0.2 cm}\\
    \vdots &q_{1} & q_{0} & \ddots & \vdots \vspace{0.2 cm}\\
   q_{n-2} & \cdots & \ddots & \ddots & q_{-1} \vspace{0.2 cm}\\
   q_{n-1}& q_{n-2}& \cdots &q_{1}& q_{0} \vspace{0.2 cm}\\
  \end{array}
\right),
$$
i.e., {$q_{ij}=q_{i-j}$}. Assume that the diagonals
are the
Fourier coefficients of {the} function $f$, i.e.,
$${q_k}=\frac{1}{2\pi}\int_{-\pi}^{\pi}f(x)
\mathrm{e}^{-\mathrm{i}kx}\mathrm{d}x,$$
then function $f(x)$ is called the generating function of
{$\mathbf{Q}_n$}.
\end{definition}

\begin{lemma}\label{Lem.3.2}(Grenander-Szeg\"{o} Theorem
\cite{Chan2011}) For the above
Toeplitz matrix ${\mathbf{Q}_n}$,
 let $f(x)$ be a $2\pi$-periodic continuous real-valued
function defined on $\left[-\pi, \pi\right]$. Denote
$\lambda_{\min}({\mathbf{Q}_n})$ and $\lambda_{\max}({\mathbf{Q}_n})$
as the smallest and largest eigenvalues of ${\mathbf{Q}_n}$,
respectively. Then one has
$$f_{\min}\leq \lambda_{\min}({\mathbf{Q}_n})
\leq \lambda_{\max}({\mathbf{Q}_n}) \leq
f_{\max},$$ where $f_{\min}$, $f_{\max}$ are
the minimum and maximum
values of $f(x)$ on $\left[-\pi, \pi\right]$.
Moreover, if $f_{\min} < f_{\max}$, then all
eigenvalues of ${\mathbf{Q}_n}$ satisfy
$$f_{\min}< \lambda({\mathbf{Q}_n}) < f_{\max},
$$ for all $n>0$.
\end{lemma}

\begin{lemma}\label{Lem.3.3}The matrix $\mathbf{K}_\gamma$ has
 the following properties:\\
$\mathrm{(1)}$~$\mathbf{K}_\gamma$ is symmetric;\;\\
$\mathrm{(2)}$~$\mathbf{K}_\gamma$ is positive definite,
 i.e., $\mathbf{W}^{\mathsf{T}}\mathbf{K}_\gamma \mathbf{W}>0$ for
any nonzero $\mathbf{W}\in \mathds{R}^{M-1}$.
\end{lemma}

\begin{proof}
(1)~From the equation (\ref{eq.2.3}), we can easily find that $g_m^{(\gamma)}=g_{-m}^{(\gamma)}$,
which shows that the
 matrix $K_\gamma$ is symmetric.

(2)~Based on the Definition
\ref{Def:1} and the equation (\ref{eq.2.2}),
we can know the generating function of the matrix
${\mathbf{K}}_\gamma$ is
\begin{equation*}
\begin{aligned}
 g(\gamma,z)&=
\sum\limits_{m=-\infty}^{+\infty}
g_{m}^{(\gamma)}e^{\mathrm{i}mz}\\
&=
\left[1+\frac{\gamma}{24}\left(2-e^{\mathrm{i}z}
-e^{-\mathrm{i}z}\right)+\frac{\gamma\left(5\gamma+22\right)}{5760}
 \left(2-e^{\mathrm{i}z}-e^{-\mathrm{i}z}\right)^2\right]
\left(2-e^{\mathrm{i}z}-e^{-\mathrm{i}z}\right)^{\frac{\gamma}{2}}\\
&=\left[1+\frac{\gamma}{6}\sin^2\left(\frac{z}{2}\right)
+\frac{\gamma(5\gamma+22)}{360}\sin^4\left(\frac{z}{2}\right)
\right]
\left[4\sin^2\left(\frac{z}{2}\right)\right]^{\frac{\gamma}{2}}\geq0
\end{aligned}
\end{equation*}
for $z\in[-\pi,\pi]$ and $\gamma\in(1,2]$.
This implies that matrix $\mathbf{K}_\gamma$ is positive definite for $1<\gamma\leq2$
based on Lemma \ref{Lem.3.2}.
\end{proof}

\begin{lemma}\label{Lem.3.4}
There exist real symmetric positive matrix
${\mathbf{L}}^{(1)}$ such that the following inner product equality holds:
\begin{equation*}
\begin{aligned}
\left({\mathbf{A}}^{(1)}\mathbf{U},\mathbf{V}\right)=\left({\mathbf{L}}^{(1)}
\mathbf{U},{\mathbf{L}}^{(1)}\mathbf{V}\right).
\end{aligned}
\end{equation*}
\end{lemma}

\begin{proof}
From the Lemma \ref{Lem.3.3}, we know that the matrix $\mathbf{A}^{(1)}$
is a real positive definite symmetric matrix.
Hence, based on the basic properties of matrix,
there is a real orthogonal matrix $\mathbf{P}$ and a real diagonal matrix
$\mathbf{\Lambda}=\mathrm{diag}\left(\lambda\right)$, such that
\begin{equation}\label{eq.3.2}
\begin{aligned}
{\mathbf{A}}^{(1)}=\mathbf{P}\Lambda\mathbf{P}^{\mathsf{T}}
=\left(\mathbf{P}\Lambda^{\frac{1}{2}}\mathbf{P}^{\mathsf{T}}\right)^{\mathsf{T}}
\left(\mathbf{P}\Lambda^{\frac{1}{2}}\mathbf{P}^{\mathsf{T}}\right)
=\left(\mathbf{L}^{(1)}\right)^{\mathsf{T}}\mathbf{L}^{(1)},
\end{aligned}
\end{equation}
where $\Lambda^{\frac{1}{2}}=\mathrm{diag}\left(\lambda^{\frac{1}{2}}\right)$
and $\mathbf{L}^{(1)}=\mathbf{P}\Lambda^{\frac{1}{2}}\mathbf{P}^{\mathsf{T}}$.
Furthermore, with the help of the decomposition (\ref{eq.3.2}), one has
\begin{equation*}
\begin{aligned}
\left({\mathbf{A}}^{(1)}\mathbf{U},\mathbf{V}\right)=
\mathbf{V}^{\mathsf{T}}{\mathbf{A}}^{(1)}\mathbf{U}
=\mathbf{V}^{\mathsf{T}}\left(\mathbf{L}^{(1)}\right)^{\mathsf{T}}
\mathbf{L}^{(1)}\mathbf{U}
=\left(\mathbf{L}^{(1)}\mathbf{V}\right)^{\mathsf{T}}
\mathbf{L}^{(1)}\mathbf{U}
=
\left({\mathbf{L}}^{(1)}
\mathbf{U},{\mathbf{L}}^{(1)}\mathbf{V}\right).
\end{aligned}
\end{equation*}
This completes the proof.
\end{proof}

\begin{lemma}\label{Lem.3.5}(Brouwder fixed point theorem \cite{Akrivis1993})
Let $\left(H,\langle\cdot,\cdot\rangle\right)$ be a finite-dimensional inner product space,
$\left\|\cdot\right\|$ be the associated norm, and $\mathcal{H}:H\rightarrow H$ be continuous.
 Assume, moreover, that
 \begin{equation*}
\begin{aligned}
\exists \sigma>0, \;\forall z\in H, \;\left\|z\right\|=\sigma,\;
\Re\left\langle \mathcal{H}(z),z\right\rangle>0.
\end{aligned}
\end{equation*}
Then, there exists a $z^*\in H$ such that $\mathcal{G}(z^*)=0$ and
 $\left\|z^*\right\|\leq\sigma$.
\end{lemma}

\begin{lemma}\label{Lem.3.6}
Let ${\mathbf{A}}^{(1)}$ be a real $(M-1)\times(M-1)$ matrix and
satisfy the conditions of Lemma \ref{Lem.3.3}, then we have
\begin{equation*}
\begin{aligned}
\displaystyle \left\|\left(
\mathbf{I}+\frac{\tau}{2}\mathbf{A}^{(1)}
\right)^{-1}\left(
\mathbf{I}-\frac{\tau}{2}\mathbf{A}^{(1)}
\right)\right\|_\infty< 1
\end{aligned}
\end{equation*}
for any $\tau>0$.
\end{lemma}

\begin{proof}
Let $\lambda_j \left(\mathbf{A}^{(1)}\right), j = 1, 2, \cdots, M -1$ be the
eigenvalues of matrix $\mathbf{A}^{(1)}$
and
denote $\mathbf{B}=\left(
\mathbf{I}+\frac{\tau}{2}\mathbf{A}^{(1)}
\right)^{-1}\left(
\mathbf{I}-\frac{\tau}{2}\mathbf{A}^{(1)}
\right)$, then the eigenvalues of the matrix $\mathbf{B}$ are
\begin{equation*}
\begin{aligned}
\displaystyle
\lambda_j \left(\mathbf{B}\right)=\frac{1-\frac{\tau}{2}\lambda_j \left(\mathbf{A}^{(1)}\right)}
{1+\frac{\tau}{2}\lambda_j \left(\mathbf{A}^{(1)}\right)},\;j = 1, 2, \cdots, M -1.
\end{aligned}
\end{equation*}

The matrix ${\mathbf{A}}^{(1)}$
satisfies the conditions of Lemma \ref{Lem.3.3}, which imply that
$\lambda_j \left(\mathbf{A}^{(1)}\right)>0$, and
it is not hard to check that $\left|\lambda_j \left(\mathbf{B}\right)\right|<1$
for any $\tau>0$.
Therefore, this reflects that the spectral radius of the matrix $\mathbf{B}$ is less than 1,
so the result to be proved is valid.
\end{proof}

\begin{lemma}\label{Lem.3.7}
Let
$ g(x)=x+\frac{\tau}{2}f(x)$ with $f(x)=x-x^3$.
Then there holds that
\begin{equation*}
\begin{aligned}
\displaystyle \left|g(x)\right|\leq1\;for\; all\; x\in[-1,1] \;and \;\tau\in(0,1].
\end{aligned}
\end{equation*}
\end{lemma}
\begin{proof}
Due to  $g'(x)=1+\frac{\tau}{2}f'(x)\geq0$ for all $x\in[-1,1]$ and $\tau\in(0,1]$,
 then we know that $g(-1)\leq g(x)\leq g(1)$. Moreover, note that $g(-1)=-1$ and
$g(1)=1$, then we deduce that
\begin{equation*}
\begin{aligned}
-1\leq g(x)\leq 1 \;for \;all\; x\in[-1,1],
\end{aligned}
\end{equation*}
which concludes the proof.
\end{proof}

%

\begin{lemma}\label{Lem.3.9}
Let $\lambda_j\left({\mathbf{A}}^{(1)}\right),\;j=1,2,\cdots,M-1$, be the eigenvalue
of the matrix ${\mathbf{A}}^{(1)}$.
Then there holds that
\begin{align*}
\lambda_j\left({\mathbf{A}}^{(1)}\right)\leq
\left\{
\begin{aligned}
  &\frac{2\varepsilon^2}{h^\gamma} g_{0}^{(\gamma)},\;\gamma\in(1,\gamma^{*}];\\
 & \frac{2\varepsilon^2}{h^\gamma}\left(g_{0}^{(\gamma)}+2g_{2}^{(\gamma)}\right),\;
  \gamma\in(\gamma^{*},2],
\end{aligned}
\right.\;\;j=1,2,\cdots,M-1.
\end{align*}
\end{lemma}

\begin{proof}
From the Gerschgorin’s circle theorem \cite{Kincaid1991}, we know that
\begin{equation*}
\begin{aligned}
\displaystyle
\left|\lambda_j\left({\mathbf{A}}^{(1)}\right)
-\frac{\varepsilon^2}{h^\gamma}g_0^{(\gamma)}\right|\leq
\frac{\varepsilon^2}{h^\gamma}\sum_{\substack{m=-M+j\\m\neq0}}^{j-1}\left|g_m^{(\gamma)}\right|
=\frac{\varepsilon^2}{h^\gamma}\sum_{m=-M+j}^{j-1}\left|g_m^{(\gamma)}\right|
-\frac{\varepsilon^2}{h^\gamma}g_0^{(\gamma)},
\;j=1,2\cdots,M-1,
\end{aligned}
\end{equation*}
that is
\begin{equation}\label{eq.3.3}
\begin{aligned}
\displaystyle
\lambda_j\left({\mathbf{A}}^{(1)}\right)
\leq\frac{\varepsilon^2}{h^\gamma}\sum_{m=-M+j}^{j-1}\left|g_m^{(\gamma)}\right|\leq
\frac{\varepsilon^2}{h^\gamma}\sum_{m=-\infty}^{+\infty}\left|g_m^{(\gamma)}\right|
:=\frac{\varepsilon^2}{h^\gamma}S_m.
\end{aligned}
\end{equation}

It follows from the Lemma \ref{Lem.2.2}, one has
\begin{equation}\label{eq.3.4}
\begin{aligned}
\displaystyle
S_m
=\sum_{m=-\infty}^{-1}\left|g_m^{(\gamma)}\right|
+\left|g_{0}^{(\gamma)}\right|
+\sum_{m=1}^{+\infty}\left|g_m^{(\gamma)}\right|
=-\sum_{m=-\infty}^{-1}g_m^{(\gamma)}
+g_{0}^{(\gamma)}
-\sum_{m=1}^{+\infty}g_m^{(\gamma)}=2g_{0}^{(\gamma)}
\end{aligned}
\end{equation}
for $\gamma\in(1,\gamma^{*}]$, and
\begin{equation}\label{eq.3.5}
\begin{aligned}
\displaystyle
S_m
=-\sum_{m=-\infty}^{-3}g_m^{(\gamma)}
+g_{-2}^{(\gamma)}-g_{-1}^{(\gamma)}+g_{0}^{(\gamma)}
-g_{1}^{(\gamma)}+g_{2}^{(\gamma)}
-\sum_{m=3}^{+\infty}g_m^{(\gamma)}
=2\left(g_{0}^{(\gamma)}+2g_{2}^{(\gamma)}\right)
\end{aligned}
\end{equation}
for $\gamma\in(\gamma^{*},2]$.
Combining (\ref{eq.3.3}), (\ref{eq.3.4}) and (\ref{eq.3.5}), we can easily know that the
conclusion to be proved is valid.
 This completes the proof.
\end{proof}

\begin{lemma}\label{Lem.3.10}(Gronwall inequality\cite{Holte2009}).
Suppose $\left\{E^k\right\}_{k=0}^\infty$ is a nonnegative sequence
and satisfies
\begin{equation*}
E^{k+1}\leq\left(1+c\tau\right)E^k+\tau g,\;\;k=0,1,\ldots.
\end{equation*}
Then we have
\begin{equation*}
E^{k}\leq\exp\left(ck\tau\right)\left(E^0+\frac{g}{c}\right),\;\;k=1,2,\ldots,
\end{equation*}
where $c$ and $g$ are nonnegative constants.
\end{lemma}

\subsection{Unique solvability}
In this section, the uniquely solvability of the finite difference scheme
(\ref{eq.3.1}) will be given by the following theorems.
\begin{theorem}\label{Th.3.10}
The solution of the finite difference scheme (\ref{eq.3.1}) exists.
\end{theorem}
\begin{proof}
Denote $\Theta=\mathbf{U}^{k+\frac{1}{2}}$.
For a fixed $k$, the difference scheme (\ref{eq.3.1}) can be rewritten as
\begin{equation}\label{eq.3.6}
\begin{aligned}\displaystyle
\Theta=&\frac{2}{2-\tau}\left\{\mathbf{U}^{k}-\frac{\tau}{2}\mathbf{A}^{(1)}\Theta
-\frac{\tau}{4}\left[\left(\mathbf{U}^{k}\right)^3+\left(2\Theta-\mathbf{U}^{k}\right)^3\right]
\right.\\
&+\frac{\tau^2}{4}\mathbf{A}^{(1)}\left(\Theta-\mathbf{U}^{k}\right)+\left.
\frac{\tau^2}{8}\mathbf{A}^{(1)}
\left[\left(\mathbf{U}^{k}\right)^3-\left(2\Theta-\mathbf{U}^{k}\right)^3\right]
\right\}.
\end{aligned}
\end{equation}

Based on the equation (\ref{eq.3.6}),
let us define a mapping $\mathcal{H} : \mathring{\Psi}_h\rightarrow \mathring{\Psi}_h$
as follows
\begin{equation*}
\begin{aligned}\displaystyle
\mathcal{H} \left(\Theta\right)=&\Theta-\frac{2}{2-\tau}
\left\{\mathbf{U}^{k}-\frac{\tau}{2}\mathbf{A}^{(1)}\Theta
-\frac{\tau}{4}\left[\left(\mathbf{U}^{k}\right)^3+\left(2\Theta-\mathbf{U}^{k}\right)^3\right]
\right.\\
&+\frac{\tau^2}{4}\mathbf{A}^{(1)}\left(\Theta-\mathbf{U}^{k}\right)+\left.
\frac{\tau^2}{8}\mathbf{A}^{(1)}
\left[\left(\mathbf{U}^{k}\right)^3-\left(2\Theta-\mathbf{U}^{k}\right)^3\right]
\right\},
\end{aligned}
\end{equation*}
which is obviously continuous. Taking the inner product of above equation with $\Theta$
 and using Lemma \ref{Lem.3.4},  there is
\begin{equation*}
\begin{aligned}\displaystyle
\left\langle\mathcal{H} \left(\Theta\right),\Theta\right\rangle
=&\left\langle\Theta,\Theta\right\rangle-\frac{2}{2-\tau}
\left\{\left\langle\mathbf{U}^{k},\Theta\right\rangle-
\frac{\tau}{2}\left\langle\mathbf{A}^{(1)}\Theta,\Theta\right\rangle
-\frac{\tau}{4}\left\langle\left(\mathbf{U}^{k}\right)^3
+\left(2\Theta-\mathbf{U}^{k}\right)^3
,\Theta\right\rangle
\right.\\
&+\frac{\tau^2}{4}\left\langle\mathbf{A}^{(1)}\left(\Theta-\mathbf{U}^{k}\right),\Theta\right\rangle+\left.
\frac{\tau^2}{8}\left\langle\mathbf{A}^{(1)}
\left[\left(\mathbf{U}^{k}\right)^3-\left(2\Theta-\mathbf{U}^{k}\right)^3\right],\Theta\right\rangle
\right\}\\
=&\left\|\Theta\right\|^2-\frac{2}{2-\tau}
\left\langle\mathbf{U}^{k},\Theta\right\rangle+\frac{\tau}{2-\tau}\left\|\mathbf{L}^{(1)}\Theta\right\|^2
+\frac{4\tau}{2-\tau}\left\|\Theta\right\|^4
-\frac{6\tau}{2-\tau}\left\langle\Theta^2\mathbf{U}^{k},\Theta\right\rangle\\
&+\frac{3\tau}{2-\tau}\left\|\Theta\mathbf{U}^{k}\right\|^2
-\frac{\tau^2}{2(2-\tau)}\left\|\mathbf{L}^{(1)}\Theta\right\|^2
+\frac{\tau^2}{2(2-\tau)}\left\langle\mathbf{L}^{(1)}\mathbf{U}^{k},\mathbf{L}^{(1)}\Theta\right\rangle\\
&-\frac{\tau^2}{2(2-\tau)}\left\langle\mathbf{L}^{(1)}\left(\mathbf{U}^{k}\right)^3,\mathbf{L}^{(1)}
\Theta\right\rangle
+\frac{2\tau^2}{2-\tau}\left\langle\mathbf{L}^{(1)}\Theta^3,\mathbf{L}^{(1)}\Theta\right\rangle\\
&-\frac{3\tau^2}{2-\tau}\left\langle\mathbf{L}^{(1)}\Theta^2\mathbf{U}^{k},\mathbf{L}^{(1)}
\Theta\right\rangle
+\frac{3\tau^2}{2(2-\tau)}\left\langle\mathbf{L}^{(1)}
\Theta\left(\mathbf{U}^{k}\right)^2,
\mathbf{L}^{(1)}\Theta\right\rangle.
\end{aligned}
\end{equation*}

Furthermore,
by using the Cauchy-Schwarz inequality,
we can obtain
\begin{equation*}
\begin{aligned}\displaystyle
\left\langle\mathcal{H} \left(\Theta\right),\Theta\right\rangle
\geq&\left\|\Theta\right\|\left[\left\|\Theta\right\|-\frac{2}{2-\tau}
\left\|{\mathbf{U}}^{k}\right\|
-\frac{6\tau}{2-\tau}\left\|\mathbf{U}^{k}\right\|\left\|\Theta\right\|^2
-\frac{\tau^2}{2(2-\tau)}
\left\|\mathbf{L}^{(1)}\right\|^2
\left\|\mathbf{U}^{k}\right\|\right.\\
&\left.
-\frac{\tau^2}{2(2-\tau)}
\left\|\mathbf{L}^{(1)}\right\|^2\left\|\mathbf{U}^{k}\right\|^3
-\frac{3\tau^2}{2-\tau}
\left\|\mathbf{L}^{(1)}\right\|^2\left\|\mathbf{U}^{k}\right\|\left\|\Theta\right\|^2
\right]
\end{aligned}
\end{equation*}
under the condition of $0<\tau<2$.

Let
\begin{equation*}
\begin{aligned}\displaystyle
\Delta=\tilde{\tau}^2
-6\tau\left\|\mathbf{U}^{k}\right\|\left(2+\tau\lambda_{\max}\left(\mathbf{A}^{(1)}\right)\right)
\left\{2\tilde{\tau}+\left\|\mathbf{U}^{k}\right\|
\left[4+\tau^2\lambda_{\max}\left(\mathbf{A}^{(1)}\right)
\left(1+\left\|\mathbf{U}^{k}\right\|^2\right)\right]\right\}\geq0,
\end{aligned}
\end{equation*}
and
\begin{equation*}
\begin{aligned}\displaystyle
\mu=\frac{\tilde{\tau}+
\sqrt{\Delta}}{6\tau\left\|\mathbf{U}^{k}\right\|\left(2+\tau
\lambda_{\max}\left(\mathbf{A}^{(1)}\right)\right)},
\end{aligned}
\end{equation*}
where $\tilde{\tau}=2-\tau$, then we obtain
$\left\langle\mathcal{H} \left(\Theta\right),\Theta\right\rangle
\geq\left\|\Theta\right\|\geq0$, when $\left\|\Theta\right\|=\mu$.
Thus, there exists a solution $\mu^*\in \mathring{\Psi}_h$
satisfying  $\mathcal{H}\left(\mu^*\right)= 0$.
Hence, we can draw a conclusion that the difference scheme
(\ref{eq.3.1}) exists solution
by follows form the Lemma \ref{Lem.3.5}. This completes the proof.
\end{proof}

For technical needs, the uniqueness of the solution will
 be considered after the Discrete maximum principle is proved.

\subsection{Discrete maximum principle}
In this section, we will show that the difference scheme
(\ref{eq.3.1}) preserves the
discrete maximum principle.

\begin{theorem}\label{Th.3.12}
Suppose that the initial condition $u_0(x)$ satisfies
$\max_{\mathbf{x}\in\overline{\Omega}}\left|u_0(x)\right| \leq1$, if
the time step $\tau$ satisfies
\begin{equation}\label{eq.3.14}
\begin{aligned}\displaystyle
0<\tau\leq1,
\end{aligned}
\end{equation}
then the numerical
solution $\left\{\mathbf{U}^{k+1}\right\}$ generated by
the difference scheme (\ref{eq.3.1}) also satisfies the discrete
maximum bound principle, i.e.,
$$\left\|\mathbf{U}^{k+1}\right\|_\infty\leq1,\;k=0,1,\cdots,N-1.$$
\end{theorem}

\begin{proof}
Here, we use mathematical induction method to prove it.
For the case of $k=0$, it is clearly true. Suppose that
$\left\|\mathbf{U}^n\right\|_\infty\leq1$ is true for $n=1,2,\cdots,k$. Next,
we will prove that $\left\|\mathbf{U}^{k+1}\right\|_\infty\leq1$ is also true.
From scheme (\ref{eq.3.1}), we can see that
\begin{equation}\label{eq.3.15}
\begin{aligned}\displaystyle
\left\|\mathbf{U}^{k+1}-\frac{\tau}{2}
\mathbf{f}\left(\mathbf{U}^{k+1}\right)\right\|_{\infty}
=\left\|\mathbf{B}\left(\mathbf{U}^{k}+\frac{\tau}{2}
\mathbf{f}\left(\mathbf{U}^{k}\right)\right)\right\|_{\infty},
\end{aligned}
\end{equation}
where $\mathbf{B}=\left(
\mathbf{I}+\frac{\tau}{2}\mathbf{A}^{(1)}
\right)^{-1}\left(
\mathbf{I}-\frac{\tau}{2}\mathbf{A}^{(1)}
\right)$.

Obviously,
 the left-hand side of equation (\ref{eq.3.15}) becomes
\begin{equation}\label{eq.3.16}
\begin{aligned}\displaystyle
\left\|\mathbf{U}^{k+1}-\frac{\tau}{2}
\mathbf{f}\left(\mathbf{U}^{k+1}\right)\right\|_{\infty}
=\left(1-\frac{\tau}{2}\right)\left\|\mathbf{U}^{k+1}\right\|_{\infty}
+\frac{\tau}{2}\left\|\mathbf{U}^{k+1}\right\|_{\infty}^3
\end{aligned}
\end{equation}
for $0<\tau\leq2$.

Under the condition of (\ref{eq.3.14}), using
Lemmas \ref{Lem.3.6} and \ref{Lem.3.7}, and note that $\left\|\mathbf{U}^k\right\|_\infty\leq1$,
we know that the right-hand side
 of equation (\ref{eq.3.15}) satisfies
\begin{equation}\label{eq.3.17}
\begin{aligned}\displaystyle
\left\|\mathbf{B}\left(\mathbf{U}^{k}+\frac{\tau}{2}
\mathbf{f}\left(\mathbf{U}^{k}\right)\right)\right\|_{\infty}
\leq\left\|\mathbf{B}\right\|_{\infty}
\left\|\mathbf{U}^{k}+\frac{\tau}{2}
\mathbf{f}\left(\mathbf{U}^{k}\right)\right\|_{\infty}
\leq1.
\end{aligned}
\end{equation}
Combining (\ref{eq.3.15}), (\ref{eq.3.16}) and (\ref{eq.3.17}) can lead to
\begin{equation*}
\begin{aligned}\displaystyle
\left(1-\frac{\tau}{2}\right)\left\|\mathbf{U}^{k+1}\right\|_{\infty}
+\frac{\tau}{2}\left\|\mathbf{U}^{k+1}\right\|_{\infty}^3\leq1,
\end{aligned}
\end{equation*}
that is
\begin{equation*}
\begin{aligned}\displaystyle
\left(\left\|\mathbf{U}^{k+1}\right\|_{\infty}-1\right)
\left(\left\|\mathbf{U}^{k+1}\right\|_{\infty}^2+
\left\|\mathbf{U}^{k+1}\right\|_{\infty}+\frac{2}{\tau}\right)\leq0,
\end{aligned}
\end{equation*}
which implies that $\left\|\mathbf{U}^{k+1}\right\|_{\infty}\leq1$.
 This completes the proof.
\end{proof}

\begin{theorem}
The solution of the  finite difference scheme (\ref{eq.3.1}) is unique.
\end{theorem}

\begin{proof}
Let $\mathbf{U}^k$ and $\mathbf{V}^k$ be the numerical solutions of the scheme (\ref{eq.3.1}) and denote
$$\xi^k=\mathbf{U}^k-\mathbf{V}^k,\;k=0,1,\cdots,N.$$

Based on difference scheme (\ref{eq.3.1}), it is not difficult to find that $\xi^k$ satisfies the following equation
\begin{equation}\label{eq.3.7}
\begin{aligned}\displaystyle
\left[\left(1-\frac{\tau}{2}\right)\mathbf{I}+\frac{\tau}{2}\mathbf{A}^{(1)}
\right]\mathbf{\xi}^{k+1}=\left[\left(1+\frac{\tau}{2}\right)\mathbf{I}-\frac{\tau}{2}\mathbf{A}^{(1)}
\right]\mathbf{\xi}^{k}+\alpha^{k+1/2}
,\;\;k=0,1,\cdots,N-1,
\end{aligned}
\end{equation}
where
\begin{equation*}
\begin{aligned}\displaystyle
\alpha^{k+1/2}=\frac{\tau}{2}\left[\left(\mathbf{V}^{k+1}\right)^3
-\left(\mathbf{U}^{k+1}\right)^3\right]
+\frac{\tau}{2}\left(
\mathbf{I}-\frac{\tau}{2}\mathbf{A}^{(1)}
\right)\left[\left(\mathbf{V}^{k}\right)^3-\left(\mathbf{U}^{k}\right)^3\right].
\end{aligned}
\end{equation*}

Computing the discrete inner product of (\ref{eq.3.7}) with $2\xi^{k+1/2}$ yield
\begin{equation}\label{eq.3.8}
\begin{aligned}\displaystyle
&\left\langle\left[\left(1-\frac{\tau}{2}\right)\mathbf{I}+\frac{\tau}{2}\mathbf{A}^{(1)}
\right]\mathbf{\xi}^{k+1},2\xi^{k+1/2}\right\rangle
-\left\langle\left[\left(1+\frac{\tau}{2}\right)\mathbf{I}-\frac{\tau}{2}\mathbf{A}^{(1)}
\right]\mathbf{\xi}^{k},2\xi^{k+1/2}\right\rangle\\
=&\left\langle\alpha^{k+1/2},2\xi^{k+1/2}\right\rangle.
\end{aligned}
\end{equation}

Next, we will proceed to analyze each term in equation (\ref{eq.3.8}). Firstly, for the left-hand
 term of equation (\ref{eq.3.8}), with the help of  Lemma \ref{Lem.3.4} and
 Cauchy-Schwarz inequality,
 we can see that
\begin{equation}\label{eq.3.9}
\begin{aligned}\displaystyle
&\left\langle\left[\left(1-\frac{\tau}{2}\right)\mathbf{I}+\frac{\tau}{2}\mathbf{A}^{(1)}
\right]\mathbf{\xi}^{k+1},2\xi^{k+1/2}\right\rangle
-\left\langle\left[\left(1+\frac{\tau}{2}\right)\mathbf{I}-\frac{\tau}{2}\mathbf{A}^{(1)}
\right]\mathbf{\xi}^{k},2\xi^{k+1/2}\right\rangle\\
\geq&\left(1-\tau\right)
\left\|\mathbf{\xi}^{k+1}\right\|^2
-\left(1+\tau\right)
\left\|\mathbf{\xi}^{k}\right\|^2.
\end{aligned}
\end{equation}

For the right-hand term of equation (\ref{eq.3.8}),
it holds that
\begin{equation}\label{eq.3.10}
\begin{aligned}\displaystyle
\left\langle\alpha^{k+1/2},2\xi^{k+1/2}\right\rangle
=&\frac{\tau}{2}\left\langle\left[\left(\mathbf{V}^{k+1}\right)^3-\left(\mathbf{U}^{k+1}\right)^3
\right],2\xi^{k+1/2}\right\rangle\\
&+\frac{\tau}{2}\left\langle\left(
\mathbf{I}-\frac{\tau}{2}\mathbf{A}^{(1)}
\right)\left[\left(\mathbf{V}^{k}\right)^3-\left(\mathbf{U}^{k}\right)^3\right],
2\xi^{k+1/2}\right\rangle\\
=&\Xi_1+\Xi_2.
\end{aligned}
\end{equation}

Next, we will analyze $\Xi_1$ and $\Xi_2$. For $\Xi_1$, it follows from
the Cauchy-Schwarz inequality,
we can know that
\begin{equation}\label{eq.3.11}
\begin{aligned}\displaystyle
\Xi_1
=&\frac{\tau}{2}\left\langle\left[\left(\mathbf{V}^{k+1}\right)^3-\left(\mathbf{U}^{k+1}\right)^3
\right],2\xi^{k+1/2}\right\rangle\leq\frac{3\tau}{4}
\hat{C}_{\max}^2\left(3\left\|\xi^{k+1}\right\|^2+\left\|\xi^{k}\right\|^2\right),
\end{aligned}
\end{equation}
where $\hat{C}_{\max}=\max\left\{\left\|\mathbf{U}^{k+1}\right\|,\left\|\mathbf{V}^{k+1}\right\|\right\}$.
Similarly, for $\Xi_2$, we also have
\begin{equation}\label{eq.3.12}
\begin{aligned}\displaystyle
\Xi_2\leq\frac{3\tau}{8}
\tilde{C}_{\max}^2
\left(2+\tau\left\|\mathbf{A}^{(1)}\right\|\right)
\left(3\left\|\xi^{k}\right\|^2+\left\|\xi^{k+1}\right\|^2\right),
\end{aligned}
\end{equation}
where $\tilde{C}_{\max}=\max\left\{\left\|\mathbf{U}^{k}\right\|,
\left\|\mathbf{V}^{k}\right\|\right\}$.

Combining (\ref{eq.3.10}), (\ref{eq.3.11}) and (\ref{eq.3.12}), we can conclude that
\begin{equation}\label{eq.3.13}
\begin{aligned}\displaystyle
\left\langle\alpha^{k+1/2},2\xi^{k+1/2}\right\rangle
\leq4\tau C_{\max}\left(\left\|\xi^{k}\right\|^2+\left\|\xi^{k+1}\right\|^2\right),
\end{aligned}
\end{equation}
where $C_{\max}=\left\{\frac{3}{4}
\hat{C}_{\max}^2,\frac{3}{8}\tilde{C}_{\max}^2\right\}$.

Substituting (\ref{eq.3.9}) and (\ref{eq.3.13}) into (\ref{eq.3.8}), and when
$\tau\left(4 C_{\max}+1\right)\leq1/2$,
 it can eventually lead to
\begin{equation*}
\begin{aligned}\displaystyle
\left\|\xi^{k}\right\|\leq3^{N/2}\left\|\xi^{0}\right\|=0,\;k=0,1,2,\cdots,N,
\end{aligned}
\end{equation*}
due to $\left\|\xi^{0}\right\|=0$,  which implies that $\mathbf{U}^k=\mathbf{V}^k$.
The proof is complete.
\end{proof}

\subsection{Discrete energy stability}
First, we define the following discrete energy function
 \begin{equation}\label{eq.3.18}
\begin{aligned}\displaystyle
{\mathcal{E}}_h\left(\mathbf{U}\right)=\left<
F\left(\mathbf{U}\right),\mathbf{1}\right>
+\frac{h}{2}\mathbf{U}^T \mathbf{A}^{(1)}\mathbf{U},
\end{aligned}
\end{equation}
for all $\mathbf{U}=\left[U_1,U_2,\cdots,U_{M-1}\right]^T\in \mathds{R}^{M-1}$ and
$\mathbf{1}=\left[1,1,\cdots,1\right]^T\in \mathds{R}^{M-1}$.
Here, we will prove the
energy stability of numerical scheme (\ref{eq.3.1}), that is, we will prove that the
energy function ${\mathcal{E}}_h\left(\mathbf{U}\right)$ decreases with time.

\begin{theorem}\label{Th.3.13}
Under the assumption of
\begin{align}\label{eq.3.19}
0<\tau\leq
\left\{
\begin{aligned}
  &\frac{4}{1+\sqrt{{1+\frac{16\varepsilon^2}{h^\gamma} g_{0}^{(\gamma)}}}}
  ,\;\;\;\gamma\in(1,\gamma^{*}];\vspace{0.2cm}\\
 & \frac{4}
 {1+\sqrt{{1+ \frac{16\varepsilon^2}{h^\gamma}\left(g_{0}^{(\gamma)}+2g_{2}^{(\gamma)}\right)}}}
,\;\;\;
  \gamma\in(\gamma^{*},2],
\end{aligned}
\right.
\end{align}
the solution generated by the difference method (\ref{eq.3.1}) satisfies
 the energy stability, i.e,.
 \begin{equation*}
\begin{aligned}\displaystyle
{\mathcal{E}}_h\left(\mathbf{U}^{k+1}\right)\leq{\mathcal{E}}_h\left(\mathbf{U}^{k}\right),\;k=0,1,\cdots,N-1.
\end{aligned}
\end{equation*}
\end{theorem}

\begin{proof}
Taking the discrete $L^2$ inner product of (\ref{eq.3.1}) with
$\left(\mathbf{U}^{k+1}-\mathbf{U}^k\right)^{\mathsf{T}}$ yields
\begin{equation}\label{eq.3.20}
\begin{aligned}\displaystyle
&\left<\mathbf{U}^{k+1}-\mathbf{U}^k,f\left(\mathbf{U}^{k+1}\right)+f\left(\mathbf{U}^{k}\right)\right>
-\frac{2}{\tau}\left\|\mathbf{U}^{k+1}-\mathbf{U}^k\right\|^2\\
&+\frac{\tau}{2} h\left(\mathbf{U}^{k+1}-\mathbf{U}^k\right)^{\mathsf{T}}\mathbf{A}^{(1)}\left(
 f\left(\mathbf{U}^{k+1}\right)- f\left(\mathbf{U}^{k}\right)\right)
=h\left(\mathbf{U}^{k+1}-\mathbf{U}^k\right)^{\mathsf{T}}
\mathbf{A}^{(1)}\left(\mathbf{U}^{k+1}+\mathbf{U}^{k}\right).
\end{aligned}
\end{equation}

Note that $\mathbf{A}^{(1)}$ is a symmetric and positive definite
 matrix by Lemma \ref{Lem.3.3}, then we know that
  \begin{equation}\label{eq.3.21}
\begin{aligned}\displaystyle
\left(\mathbf{U}^{k+1}-\mathbf{U}^k\right)^{\mathsf{T}}\mathbf{A}^{(1)}
\left(\mathbf{U}^{k+1}+\mathbf{U}^k\right)
=\left(\mathbf{U}^{k+1}\right)^{\mathsf{T}}\mathbf{A}^{(1)}\mathbf{U}^{k+1}
-\left(\mathbf{U}^{k}\right)^{\mathsf{T}}\mathbf{A}^{(1)}\mathbf{U}^{k}.
\end{aligned}
\end{equation}
%

Considering the energy function (\ref{eq.3.18}) in two adjacent time
 levers $t_k$ and $t_{k+1}$, it is clear that
 \begin{equation}\label{eq.3.22}
\begin{aligned}\displaystyle
{\mathcal{E}}_h\left(\mathbf{U}^{k+1}\right)-{\mathcal{E}}_h\left(\mathbf{U}^k\right)=&\left<
F\left(\mathbf{U^{k+1}}\right)-
F\left(\mathbf{U^{k}}\right),\mathbf{1}\right>\\&
+\frac{h}{2}\left[\left(\mathbf{U}^{k+1}\right)^{\mathsf{T}}\mathbf{A}^{(1)}
\mathbf{U}^{k+1}
-
\left(\mathbf{U}^k\right)^{\mathsf{T}}\mathbf{A}^{(1)}\mathbf{U}^k\right].
\end{aligned}
\end{equation}
Combining (\ref{eq.3.20}), (\ref{eq.3.21}) and (\ref{eq.3.22}) can lead to
 \begin{equation}\label{eq.3.23}
\begin{aligned}\displaystyle
&{\mathcal{E}}_h\left(\mathbf{U}^{k+1}\right)-{\mathcal{E}}_h\left(\mathbf{U}^k\right)\\=&\left<
F\left(\mathbf{U^{k+1}}\right)-
F\left(\mathbf{U^{k}}\right),\mathbf{1}\right>
+\frac{1}{2}\left<\mathbf{U}^{k+1}-\mathbf{U}^k,f\left(\mathbf{U}^{k+1}\right)+f\left(\mathbf{U}^{k}\right)\right>\\&
-\frac{1}{\tau}\left\|\mathbf{U}^{k+1}-\mathbf{U}^k\right\|^2+\frac{\tau}{4} h\left(\mathbf{U}^{k+1}
-\mathbf{U}^k\right)^{\mathsf{T}}\mathbf{B}^{(1)}\left(
 f\left(\mathbf{U}^{k+1}\right)- f\left(\mathbf{U}^{k}\right)\right).
\end{aligned}
\end{equation}
It follows from the following fundamental inequalities
 \begin{equation*}
\begin{aligned}\displaystyle
\frac{1}{4}\left[\left(\alpha^2-1\right)^2-\left(\beta^2-1\right)^2\right]\leq
\left(\alpha^3-\alpha\right)\left(\alpha-\beta\right)+\frac{1}{2}\left(\alpha-\beta\right)^2,
\end{aligned}
\end{equation*}
and
 \begin{equation*}
\begin{aligned}\displaystyle
\frac{1}{4}\left[\left(\alpha^2-1\right)^2-\left(\beta^2-1\right)^2\right]\leq
\left(\beta^3-\beta\right)\left(\alpha-\beta\right)+\frac{1}{2}\left(\alpha-\beta\right)^2,
\end{aligned}
\end{equation*}
for any $\alpha,\beta\in[-1,1]$, we get
 \begin{equation}\label{eq.3.24}
\begin{aligned}\displaystyle
\left<
F\left(\mathbf{U^{k+1}}\right)-
F\left(\mathbf{U^{k}}\right),\mathbf{1}\right>
=&\frac{1}{4}h\sum_{j=1}^{M-1}\left[
\left(\left(U^{k+1}_j\right)^2-
1\right)^2-\left(\left(U^{k}_j\right)^2-
1\right)^2
\right]
\\
\leq&\frac{1}{2}h\sum_{j=1}^{M-1}\left[
\left(\left(U^{k+1}_j\right)^3-
U^{k+1}_j\right)\left(U^{k+1}_j-
U^{k}_j\right)\right.\\&\left.
+\left(\left(U^{k}_j\right)^3-
U^{k}_j\right)\left(U^{k+1}_j-
U^{k}_j\right)
+\left(U^{k+1}_j-
U^{k}_j\right)^2
\right]\\
=&-\frac{1}{2}\left<\mathbf{U}^{k+1}-\mathbf{U}^k,
f\left(\mathbf{U}^{k+1}\right)+f\left(\mathbf{U}^{k}\right)\right>
+\frac{1}{2}\left\|\mathbf{U}^{k+1}-
\mathbf{U}^{k}\right\|^2.
\end{aligned}
\end{equation}
Substituting (\ref{eq.3.24}) into (\ref{eq.3.23}), we further get
 \begin{equation}\label{eq.3.25}
\begin{aligned}\displaystyle
&{\mathcal{E}}_h\left(\mathbf{U}^{k+1}\right)-{\mathcal{E}}_h\left(\mathbf{U}^k\right)\\
\leq&
\left(\frac{1}{2}-\frac{1}{\tau}\right)\left\|\mathbf{U}^{k+1}-\mathbf{U}^k\right\|^2
+\frac{\tau}{4} h\left(\mathbf{U}^{k+1}
-\mathbf{U}^k\right)^{\mathsf{T}}\mathbf{B}^{(1)}\left(
 f\left(\mathbf{U}^{k+1}\right)- f\left(\mathbf{U}^{k}\right)\right) \\
 \leq&
\left(\frac{1}{2}-\frac{1}{\tau}\right)\left\|\mathbf{U}^{k+1}-\mathbf{U}^k\right\|^2
+\frac{\tau}{4} \left\|\mathbf{U}^{k+1}
-\mathbf{U}^k\right\|\left\|\mathbf{B}^{(1)}\right\|\left\|
 f\left(\mathbf{U}^{k+1}\right)- f\left(\mathbf{U}^{k}\right)\right\|.
\end{aligned}
\end{equation}

Since $\left\|U^k\right\|_\infty \leq1$  for $k= 0, 1, \cdots , N$, by the mean value theorem,
we know that
there exist
$\xi_1, \xi_2, \cdots , \xi_{M-1}$, such that $|\xi_j| \leq1$ for $j= 1, 2,\cdots , M-1$ and
 \begin{equation}\label{eq.3.26}
\begin{aligned}\displaystyle
f\left({U_j^{k+1}}\right)-
f\left({U_j^{k}}\right)=f'\left(\xi_j\right)\left(U_j^{k+1}-U_j^{k}\right),
j= 1, 2,\cdots , M-1.
\end{aligned}
\end{equation}
Denoting $\left\|f'\right\|_{C[-1,1]}=\max_{\xi_j\in[-1,1]}\left|f'\left(\xi_j\right)\right|$.
By applying the equation (\ref{eq.3.26}), we can see that
 \begin{equation}\label{eq.3.27}
\begin{aligned}\displaystyle
\left\|f\left(\mathbf{U}^{k+1}\right)- f\left(\mathbf{U}^{k}\right)\right\|=&
\sqrt{h\sum_{j=1}^{M-1}\left(f\left({U}^{k+1}_j\right)- f\left({U}^{k}_j\right)\right)^2}
=\sqrt{h\sum_{j=1}^{M-1}f'^2\left(\xi_j\right)\left(U_j^{k+1}-U_j^{k}\right)^2}\\
\leq&\left\|f'\right\|_{C[-1,1]}\left\|\mathbf{U}^{k+1}
-\mathbf{U}^k\right\|.
\end{aligned}
\end{equation}
Substituting (\ref{eq.3.27}) into (\ref{eq.3.25}) can result in
 \begin{equation*}
\begin{aligned}\displaystyle
{\mathcal{E}}_h\left(\mathbf{U}^{k+1}\right)-{\mathcal{E}}_h\left(\mathbf{U}^k\right)
 \leq&
\left(\frac{1}{2}-\frac{1}{\tau}+\frac{\tau}{4}\left\|\mathbf{A}^{(1)}\right\|\left\|f'\right\|_{C[-1,1]}
\right)\left\|\mathbf{U}^{k+1}-\mathbf{U}^k\right\|^2\\
=&\left(\frac{1}{2}-\frac{1}{\tau}+\frac{\tau}{4}\lambda_{\max}\left(\mathbf{A}^{(1)}\right)\left\|f'\right\|_{C[-1,1]}
\right)\left\|\mathbf{U}^{k+1}-\mathbf{U}^k\right\|^2.
\end{aligned}
\end{equation*}
Therefore, we can get the final result
\begin{equation*}
\begin{aligned}\displaystyle
{\mathcal{E}}_h\left(\mathbf{U}^{k+1}\right)
\leq{\mathcal{E}}_h\left(\mathbf{U}^k\right),\;k=0,1,2\cdots,N-1,
\end{aligned}
\end{equation*}
with the help of the Lemma \ref{Lem.3.9} and condition (\ref{eq.3.19}).
This completes the proof.
\end{proof}

\subsection{Convergence analysis}
In this section, we will give the error estimation of the fully discrete scheme (\ref{eq.3.1})
in the discrete $L^\infty$ norm. Defining the following error function
\begin{equation*}
\begin{aligned}\displaystyle
{\mathbf{E}}^{k}=\mathbf{u}(x,t_k)-
\mathbf{U}^k,\;k=0,1,2\cdots,N.
\end{aligned}
\end{equation*}

From the process of establishing the difference scheme earlier, it is not difficult to find that
the truncation error of the scheme (\ref{eq.3.1}) is $R=\left(R_h+R_{\tau}\right)\tau$,
and which given by
\begin{equation*}
\begin{aligned}\displaystyle
R=\mathbf{u}\left(x,t_{k+1}\right)-\mathbf{B}\mathbf{u}\left(x,t_{k}\right)
-\frac{\tau}{2}\mathbf{g}\left(\mathbf{u}\left(x,t_{k+1}\right)\right)-
\frac{\tau}{2}\mathbf{B}
\mathbf{g}\left(\mathbf{u}\left(x,t_{k}\right)\right).
\end{aligned}
\end{equation*}

Moreover, if the exact solution is smooth enough and belongs to
$C^2\left([0,T];\mathscr{C}^{6+\gamma}\left(\overline{\Omega}\right)\right)$,
then there holds that
\begin{equation*}
\begin{aligned}\displaystyle
\left\|R\right\|_\infty\leq\left(\left\|R_\tau\right\|_\infty+\left\|R_h\right\|_\infty\right)\tau
\leq\max\left\{2\tau C_3, 2\tau C_4\right\}\left(\tau^2+h^6\right).
\end{aligned}
\end{equation*}
by follows from the (\ref{eq.2.12}) and (\ref{eq.2.15}).

\begin{theorem}\label{Th.3.14}
If the exact solution of the problem (\ref{eq.3.1}) belongs to
$C^2\left([0,T];\mathscr{C}^{6+\gamma}\left(\overline{\Omega}\right)\right)$
and initial value $u_0(x)$ meets
  $\max_{{x}\in\overline{\Omega}}\left|u_0({x})\right| \leq1$.
  Then, we have the following error estimate for the
  scheme (\ref{eq.3.1}) as
\begin{equation*}
\begin{aligned}\displaystyle
\left\|{\mathbf{E}}^{k}\right\|_\infty\leq4e^T\max\left\{C_3,C_4\right\}
\left(\tau^2+h^6\right),\;k=0,1,2\cdots,N,
\end{aligned}
\end{equation*}
under the condition (\ref{eq.3.14}).
\end{theorem}

\begin{proof}
Subtracting (\ref{eq.2.16}) from (\ref{eq.2.14}) and taking $d=1$, then we
can get the following error equation
\begin{equation}\label{eq.3.28}
\begin{aligned}\displaystyle
\mathbf{E}^{k+1}=&\left(
\mathbf{I}+\frac{\tau}{2}\mathbf{A}^{(1)}
\right)^{-1}\left(
\mathbf{I}-\frac{\tau}{2}\mathbf{A}^{(1)}
\right)\mathbf{E}^{k}+
\frac{\tau}{2}
\left[\mathbf{f}\left(\mathbf{u}\left(x,t_{k+1}\right)\right)
-\mathbf{f}\left(\mathbf{U}^{k+1}\right)\right]\\&+
\frac{\tau}{2}\left(
\mathbf{I}+\frac{\tau}{2}\mathbf{A}^{(1)}
\right)^{-1}\left(
\mathbf{I}-\frac{\tau}{2}\mathbf{A}^{(1)}
\right)\left[\mathbf{f}\left(\mathbf{u}\left(x,t_{k}\right)\right)
-\mathbf{f}\left(\mathbf{U}^{k}\right)\right]+R
.
\end{aligned}
\end{equation}
Reorganizing equation (\ref{eq.3.28}), it is easy to obtain
\begin{equation}\label{eq.3.29}
\begin{aligned}\displaystyle
&\left[\left(1-\frac{\tau}{2}\right)\mathbf{E}^{k+1}
+\frac{\tau}{2}\left(\left(\mathbf{U}^{k+1}\right)^2+
\mathbf{U}^{k+1}\mathbf{u}\left(x,t_{k+1}\right)
+\mathbf{u}^2\left(x,t_{k+1}\right)\right)\mathbf{E}^{k+1}
\right]\\
=&\mathbf{B}
\left[\left(1+\frac{\tau}{2}\right)\mathbf{E}^{k}
-\frac{\tau}{2}\left(\left(\mathbf{U}^{k}\right)^2+
\mathbf{U}^{k}\mathbf{u}\left(x,t_{k}\right)
+\mathbf{u}^2\left(x,t_{k}\right)\right)\mathbf{E}^{k}
+R\right],
\end{aligned}
\end{equation}
where $\mathbf{B}=\left(
\mathbf{I}+\frac{\tau}{2}\mathbf{A}^{(1)}
\right)^{-1}\left(
\mathbf{I}-\frac{\tau}{2}\mathbf{A}^{(1)}
\right)$.

Considering that each element of  $\left(1-\frac{\tau}{2}\right)\mathbf{I}
+\frac{\tau}{2}\left(\left(\mathbf{U}^{k+1}\right)^2+
\mathbf{U}^{k+1}\mathbf{u}\left(x,t_{k+1}\right)
+\mathbf{u}^2\left(x,t_{k+1}\right)\right)$ is of the form
$
g(x,y)=\left(1-\frac{\tau}{2}\right)
+\frac{\tau}{2}\left(x^2+
xy
+y^2\right),
$
and a simple computation shows that
\begin{equation}\label{eq.3.30}
\begin{aligned}\displaystyle
g_1(x,y)=\left(1-\frac{\tau}{2}\right)
+\frac{\tau}{2}\left(x^2+
xy
+y^2\right)\geq1-\frac{\tau}{2}
\end{aligned}
\end{equation}
for $x,y\in[-1,1]$ and $0<\tau\leq2$.
Hence, as $\left\|\mathbf{U}^{k+1}\right\|_\infty\leq1$
and $\left\|u\left(x,t_{k+1}\right)\right\|_\infty\leq1$,
it follows from (\ref{eq.3.30}) that
the left-hand side of equation (\ref{eq.3.29}) can be further reduced to
\begin{equation}\label{eq.3.31}
\begin{aligned}\displaystyle
\left\|\mathrm{LHS}\right\|_\infty
\geq\left(1-\frac{\tau}{2}\right)\left\|\mathbf{E}^{k+1}\right\|_\infty.
\end{aligned}
\end{equation}

Similarly, we can also know that
\begin{equation}\label{eq.3.32}
\begin{aligned}\displaystyle
g_2(x,y)=\left(1+\frac{\tau}{2}\right)
-\frac{\tau}{2}\left(x^2+
xy
+y^2\right)\leq1+\frac{\tau}{2}
\end{aligned}
\end{equation}
when $x,y\in[-1,1]$ and $\tau>0$.

On the other hand, based on Lemma \ref{Lem.3.6} and (\ref{eq.3.32}), the right-hand
side of equation (\ref{eq.3.29})
can be estimated as
\begin{equation}\label{eq.3.33}
\begin{aligned}\displaystyle
\left\|\mathrm{RHS}\right\|_\infty=&\left\|\mathbf{B}
\left[\left(1+\frac{\tau}{2}\right)\mathbf{E}^{k}
-\frac{\tau}{2}\left(\left(\mathbf{U}^{k}\right)^2+
\mathbf{U}^{k}\mathbf{u}\left(x,t_{k}\right)
+\mathbf{u}^2\left(x,t_{k}\right)\right)\mathbf{E}^{k}+R
\right]\right\|_\infty\\
\leq&
\left\|\left[\left(1+\frac{\tau}{2}\right)\mathbf{I}
-\frac{\tau}{2}\left(\left(\mathbf{U}^{k}\right)^2+
\mathbf{U}^{k}\mathbf{u}\left(x,t_{k}\right)
+\mathbf{u}^2\left(x,t_{k}\right)\right)
\right]\mathbf{E}^{k}+R\right\|_\infty\\
\leq&\left(1+\frac{\tau}{2}\right)\left\|\mathbf{E}^{k}\right\|_\infty
+\max\left\{2\tau C_3, 2\tau C_4\right\}\left(\tau^2+h^6\right)
\end{aligned}
\end{equation}
for  $\left\|\mathbf{U}^{k}\right\|_\infty\leq1$
and $\left\|u\left(x,t_{k}\right)\right\|_\infty\leq1$.

Combining (\ref{eq.3.31}) and (\ref{eq.3.33}), and using the condition (\ref{eq.3.14}),
then it can lead to the following result
\begin{equation*}
\begin{aligned}\displaystyle
\left\|\mathbf{E}^{k+1}\right\|_\infty\leq\left(1+{\tau}\right)\left\|\mathbf{E}^{k}\right\|_\infty
+\max\left\{4\tau C_3, 4\tau C_4\right\}\left(\tau^2+h^6\right),\;k=0,1,\cdots,N-1.
\end{aligned}
\end{equation*}
Using Lemma \ref{Lem.3.9} to the above inequality again, it is not difficult to obtain
\begin{equation*}
\begin{aligned}\displaystyle
\left\|\mathbf{E}^{k}\right\|_\infty \leq e^T\left[\left\|\mathbf{E}^{0}\right\|_\infty
+4\max\left\{C_3,C_4\right\}\left(\tau^2+h^6\right)\right],\;k=0,1,\cdots,N.
\end{aligned}
\end{equation*}
Note the fact that $\left\|\mathbf{E}^{0}\right\|_\infty=0$, so the result to be
 proved is naturally obtained. This completes the proof.
\end{proof}

\section{Numerical examples}

In this section, some numerical experiments are presented to
demonstrate the theoretical results obtained in the previous and
efficiency of the proposed numerical scheme, which are
 mainly divided into two parts.
The first one
is to test the accuracy of the numerical differential formula
(\ref{eq.2.6}) with (\ref{eq.2.7}) and
the fully discrete scheme (\ref{eq.3.1}), respectively.
And
the second part is to study the discrete maximum
 principle of the numerical scheme
 (\ref{eq.3.1}) under the time-step
 restriction.

\subsection{Convergence order tests}

In this subsection, in order to test the  convergence orders for the numerical differential formula
(\ref{eq.2.6}) with (\ref{eq.2.7}) and
the fully discrete numerical scheme (\ref{eq.3.1}),
we use exact solutions with sufficient regularity as the testing examples.

\begin{example}\label{Ex.1}
Check the convergence order of the numerical differential formula
(\ref{eq.2.6}) with (\ref{eq.2.7}).

Here, we choose the test function $u(x)=x^4(1-x)^4$ for $x\in(0,1)$.
Using the numerical differential formula
(\ref{eq.2.6}) with (\ref{eq.2.7}) to treat this test function, and the results
are shown in Table \ref{Tab.1}. From these numerical results, it is not difficult to
find that the numerical differential formula
(\ref{eq.2.6}) with (\ref{eq.2.7}) has a convergence order
of 6.

\begin{table}[htbp]\renewcommand\arraystretch{1.5}
 \begin{center}
 \caption{The maximum absolute errors and convergence orders of the formula (\ref{eq.2.6})
 with (\ref{eq.2.7}) for $\gamma\in(1,2]$.}\label{Tab.1}
 \begin{footnotesize}
\begin{tabular}{c c c c c c }\hline
  $\gamma$ &$h$&\textrm{The maximum absolute errors}&\textrm {The convergence orders}\\\hline
 $1.1 $& $\frac{1}{20}$ &8.755661e-008  &  ---\\
  $$  & $\frac{1}{30}$& 	7.726314e-009	&  	5.9873\\
  $$& $\frac{1}{40}$ &	1.377568e-009 &  5.9938 \\
  $$& $\frac{1}{50}$ & 3.614164e-010	 &   5.9963\\
    $$& $\frac{1}{60}$ &1.210898e-010 &  	5.9976\\ \hline
    $1.3 $& $\frac{1}{20}$ &1.625749e-007 &  ---\\
  $$  & $\frac{1}{30}$& 	1.433460e-008	& 5.9893\\
  $$& $\frac{1}{40}$ &	2.555080e-009 &  5.9948\\
  $$& $\frac{1}{50}$ &6.702622e-010	 &  5.9969\\
    $$& $\frac{1}{60}$ &2.245534e-010	 &5.9979\\ \hline
   $1.5 $& $\frac{1}{20}$ &2.921865e-007 &  ---\\
  $$  & $\frac{1}{30}$& 2.573720e-008	&  5.9918\\
  $$& $\frac{1}{40}$ &4.585961e-009  & 5.9960 \\
  $$& $\frac{1}{50}$ &1.202822e-009&  5.9976\\
    $$& $\frac{1}{60}$ &4.029377e-010	 & 	5.9984\\ \hline
   $1.7 $& $\frac{1}{20}$ &5.106500e-007 &  ---\\
  $$  & $\frac{1}{30}$& 	4.492734e-008	&  	5.9947\\
 $$& $\frac{1}{40}$ &	8.002056e-009	 & 	5.9974 \\
$$& $\frac{1}{50}$ & 2.098410e-009	& 	5.9985\\
 $$& $\frac{1}{60}$ &7.028819e-010	 &  	5.9990\\ \hline
    $1.9 $& $\frac{1}{20}$ &5.106500e-007  &  ---\\
  $$  & $\frac{1}{30}$& 4.492734e-008 &  5.9947\\
  $$& $\frac{1}{40}$ &		8.002056e-009 & 5.9974 \\
  $$& $\frac{1}{50}$ & 2.098410e-009	 & 5.9985\\
    $$& $\frac{1}{60}$ &	7.028819e-010	&  5.9990\\ \hline
      $2 $& $\frac{1}{20}$ &1.125000e-006  &  ---\\
  $$  & $\frac{1}{30}$& 	9.876543e-008	& 6.0000\\
  $$& $\frac{1}{40}$ &1.757812e-008	 &  6.0000 \\
  $$& $\frac{1}{50}$ & 	4.608001e-009 &  6.0000\\
    $$& $\frac{1}{60}$ &1.543202e-009	 &6.0000\\ \hline
\end{tabular}
 \end{footnotesize}
 \end{center}
 \end{table}

In addition, in Table \ref{Tab.2}, we also give the results for case $\gamma\in(0,1)$,
which shown that its convergence order can also reach sixth-order.
Combining Table \ref{Tab.1} and Table \ref{Tab.2}, we can claim that the
 numerical numerical differential formula
(\ref{eq.2.6}) with (\ref{eq.2.7}) can achieve a convergence
 order of sixth-order for all $\gamma\in(0,1)\cup(1,2]$, which is superior to
 existing some numerical
 differential formulas, as they are only valid in the
 case of $\gamma\in(0,1)$ or $\gamma\in(1,2]$.

\begin{table}[htbp]\renewcommand\arraystretch{1.5}
 \begin{center}
 \caption{The maximum absolute errors and convergence orders of the formula (\ref{eq.2.6})
 with (\ref{eq.2.7}) for $\gamma\in(0,1)$.}\label{Tab.2}
 \begin{footnotesize}
\begin{tabular}{c c c c c c }\hline
  $\gamma$ &$h$&\textrm{The maximum absolute errors}&\textrm {The convergence orders}\\\hline
 $0.1 $& $\frac{1}{20}$ &7.696507e-010  &  ---\\
  $$  & $\frac{1}{30}$& 	6.801712e-011	&  	5.9837\\
  $$& $\frac{1}{40}$ &	1.213353e-011 &  5.9920 \\
  $$& $\frac{1}{50}$ & 3.184138e-012	 &  5.9952\\
    $$& $\frac{1}{60}$ &1.066997e-012 &  		5.9967\\ \hline
    $0.3 $& $\frac{1}{20}$ &3.695955e-009	 &  ---\\
  $$  & $\frac{1}{30}$&3.266415e-010	&	5.9836\\
  $$& $\frac{1}{40}$ &5.827008e-011 &	5.9919\\
  $$& $\frac{1}{50}$ &1.529151e-011 & 5.9952\\
    $$& $\frac{1}{60}$ &	5.124070e-012 &5.9968\\ \hline
   $0.5 $& $\frac{1}{20}$ &9.873718e-009 &  ---\\
  $$  & $\frac{1}{30}$&8.725086e-010	&5.9839\\
  $$& $\frac{1}{40}$ &	1.556410e-010	  & 5.9921 \\
  $$& $\frac{1}{50}$ &	4.084322e-011	& 5.9953\\
    $$& $\frac{1}{60}$ &1.368615e-011		 & 	5.9969\\ \hline
   $0.7 $& $\frac{1}{20}$ &	2.210839e-008&  ---\\
  $$  & $\frac{1}{30}$& 	1.953068e-009		&  5.9846\\
 $$& $\frac{1}{40}$ &		3.483578e-010		 &	5.9925 \\
$$& $\frac{1}{50}$ & 	9.141122e-011		& 	5.9955\\
 $$& $\frac{1}{60}$ &3.062998e-011& 5.9970\\ \hline
    $0.9 $& $\frac{1}{20}$ &4.526141e-008	  &  ---\\
  $$  & $\frac{1}{30}$& 	3.996564e-009 & 	5.9858\\
  $$& $\frac{1}{40}$ &		7.127289e-010	 &	5.9930 \\
  $$& $\frac{1}{50}$ &1.870102e-010	 &	5.9959\\
    $$& $\frac{1}{60}$ &	6.266051e-011	&  	5.9973\\ \hline
\end{tabular}
 \end{footnotesize}
 \end{center}
 \end{table}
\end{example}

\begin{example}\label{Ex.2}
Test the convergence order of the numerical scheme
(\ref{eq.3.1}). For this purpose, we consider
the following one-dimensional space fractional
Allen-Cahn equation with a force term, i.e.
\begin{equation*}
\begin{aligned}\displaystyle
\partial_tu=\epsilon^2\mathscr{L}^\gamma u
+u-u^3+s(x,t),\;\;(x,t)\in\Omega\times(0,1],
\end{aligned}
\end{equation*}
where $\Omega =[0,1]$ and the source term is
\begin{equation*}
\begin{aligned}\displaystyle
s(x,t)=&e^{-3 t}
x^{18}(1-x)^{18}
-2e^{-t}x^{6}(1-x)^{6}\displaystyle\vspace{0.2cm}\\&
+\frac{\epsilon^2e^{-t}}{2\cos\left(\frac{\pi}{2}\gamma\right)}
\sum_{\ell=0}^{6}\frac{(-1)^{\ell}6!\,(6+\ell)!}{\ell!\,(6-\ell)!\,\Gamma\left(7+\ell-\beta\right)}
\left[x^{6+\ell-\beta}+\left(1-x\right)^{6+\ell-\beta}\right].
\end{aligned}
\end{equation*}
The exact solution of the above problem is $u(x, t)=e^{-t}x^6(1-x)^6$.

Table \ref{Tab.3} shows the maximum absolute errors and convergence orders of various
fractional orders, i.e., $\gamma=1.2,1.4,1.6,1.8$, at $t = 1$ by applying
 the numerical scheme (\ref{eq.3.1}).
From the table, we can observe that the finite
difference scheme  (\ref{eq.3.1}) has second-order temporal
accuracy and sixth-order spatial accuracy,
which is in good agreement with the proven error estimates in Theorem \ref{Th.3.14}.

\begin{table}[htbp]\renewcommand\arraystretch{1.5}
 \begin{center}
 \caption{The maximum absolute error, temporal convergence
order and spatial convergence order of the numerical solutions with $\epsilon=0.001$
at $t= 1$.}\label{Tab.3}
 \begin{footnotesize}
\begin{tabular}{c c c c c c }\hline
  $\gamma$ &$\tau$,\;$h$&The maximum absolute error& Temporal convergence
order& Spatial convergence order\\\hline
  $1.2$& $\frac{1}{8},\frac{1}{8}$ &  8.677284e-007		&  ---&  ---\\
  $$  & $\frac{1}{64},\frac{1}{16}$&       1.372736e-008		 &   1.9940& 5.9821\\
  $$& $\frac{1}{512},\frac{1}{32}$ &     3.494463e-010 	&    1.7653&   5.2958\\ \hline
  $1.4$& $\frac{1}{8},\frac{1}{8}$ &    8.681616e-007  		&  ---&  ---\\
  $$  & $\frac{1}{64},\frac{1}{16}$&    1.377418e-008     &  1.9926&    5.9779\\
  $$& $\frac{1}{512},\frac{1}{32}$ &     3.551702e-010 			&     1.7591&   5.2773\\ \hline
 $1.6$& $\frac{1}{8},\frac{1}{8}$ &     8.688534e-007	&  ---&  ---\\
  $$  & $\frac{1}{64},\frac{1}{16}$&      1.384790e-008   		 &   1.9905&   5.9714\\
  $$& $\frac{1}{512},\frac{1}{32}$ &   3.641537e-010  		&    1.7497&   5.2490\\ \hline
 $1.8$& $\frac{1}{8},\frac{1}{8}$ &   8.699670e-007			&  ---&  ---\\
  $$  & $\frac{1}{64},\frac{1}{16}$&     1.396475e-008  		 &   1.9870&   5.9611\\
  $$& $\frac{1}{512},\frac{1}{32}$ &    3.783460e-010 	&      1.7353&     5.2059\\ \hline
\end{tabular}
 \end{footnotesize}
 \end{center}
 \end{table}

Furthermore, the evolutions of the error surface are also plotted
 in Figs. \ref{Fig.21} and \ref{Fig.22} for different $\gamma$.
From these figures, we can find that the errors are very small,
which shows that our difference scheme is effective and has high-order convergence.

\begin{figure}[!htb]
\centering
 \includegraphics[width=10 cm]{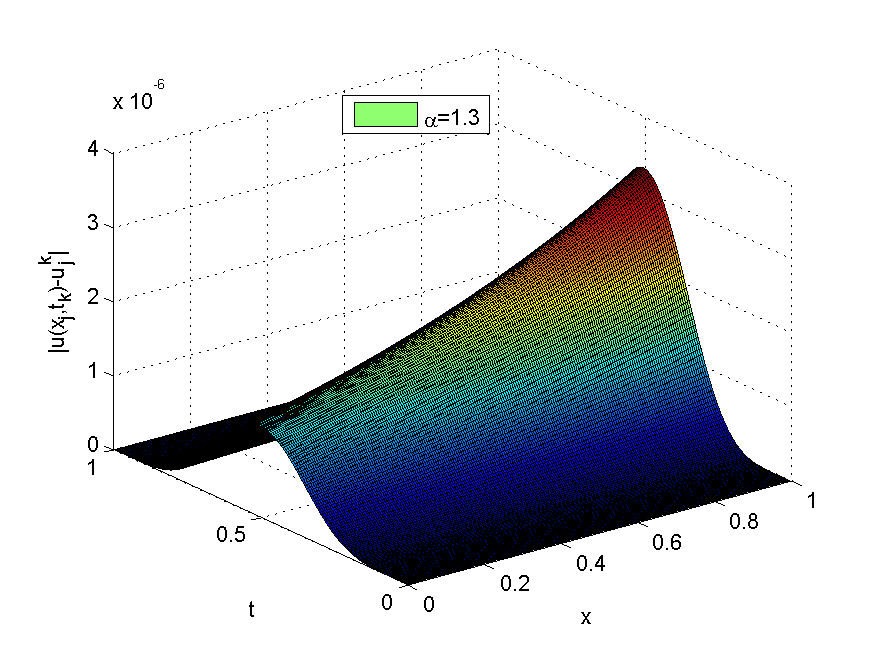}\\
  \caption{Evolutions of the error surface for
   $\gamma=1.3$ and $\epsilon=0.005$, where $\tau=0.005$ and $h=0.005$.}\label{Fig.21}
\end{figure}

\begin{figure}[!htb]
\centering
 \includegraphics[width=10 cm]{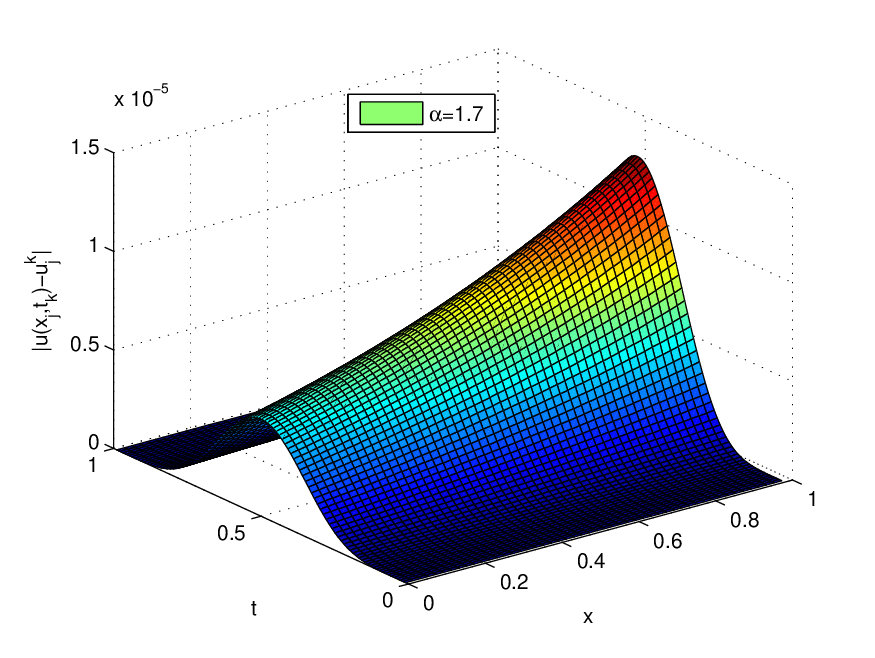}\\
  \caption{Evolutions of the error surface for
   $\gamma=1.7$ and $\epsilon=0.01$, where $\tau=0.02$ and $h=0.001$.}\label{Fig.22}
\end{figure}

\subsection{Discrete maximum principle tests}

\begin{example}\label{Ex.3}
Verify the discrete maximum principle of the numerical scheme (\ref{eq.3.1}).
Here, we consider the one-dimensional space fractional Allen-Cahn equation (\ref{eq.1.1})
with the initial condition
\begin{equation*}
\begin{aligned}\displaystyle
u_0(x)=x^{3.5+\gamma}(1-x)^{3.5+\gamma}\sin\pi x,\;\;x\in[0,1],
\end{aligned}
\end{equation*}
where the zero boundary values are set for the initial condition $u_0(x)$.
\end{example}

For this example, we fix $h=0.01$ and $\epsilon=0.1$, but vary $\tau$.
The evolution graphs of the maximum values of numerical solution
with different $\gamma$ are given Figs. \ref{Fig.31}, \ref{Fig.32}, \ref{Fig.33}
and \ref{Fig.34}, respectively.
From these figures, we can see that they are all bounded by $1$, which is consistent
with our theoretical analysis and further implies that condition (\ref{eq.3.14}) in
Theorem \ref{Th.3.12} is reasonable.

\begin{figure}[!htb]
\centering
 \includegraphics[width=10 cm]{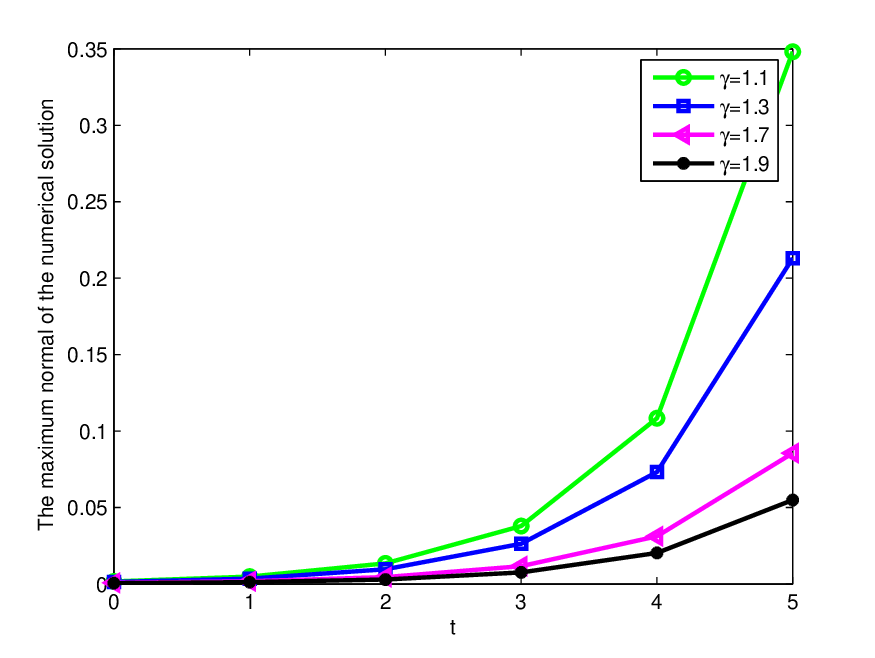}\\
  \caption{Evolutions of the maximum norm of the
  numerical solutions with different $\gamma$ values for $\tau=1$.}\label{Fig.31}
\end{figure}

\begin{figure}[!htb]
\centering
 \includegraphics[width=10 cm]{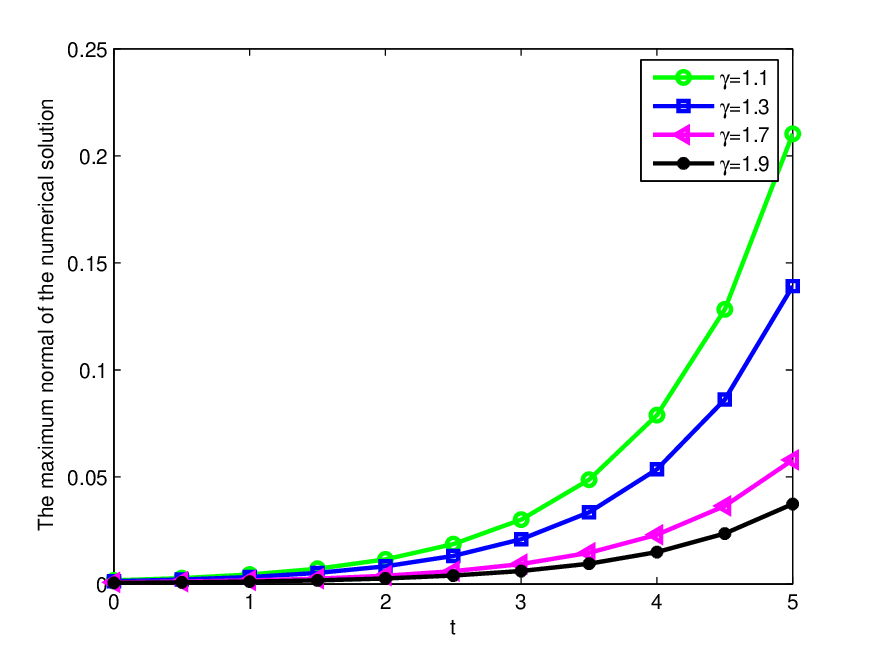}\\
  \caption{Evolutions of the maximum norm of the
  numerical solutions with different $\gamma$ values for $\tau=0.5$.}\label{Fig.32}
\end{figure}

\begin{figure}[!htb]
\centering
 \includegraphics[width=10 cm]{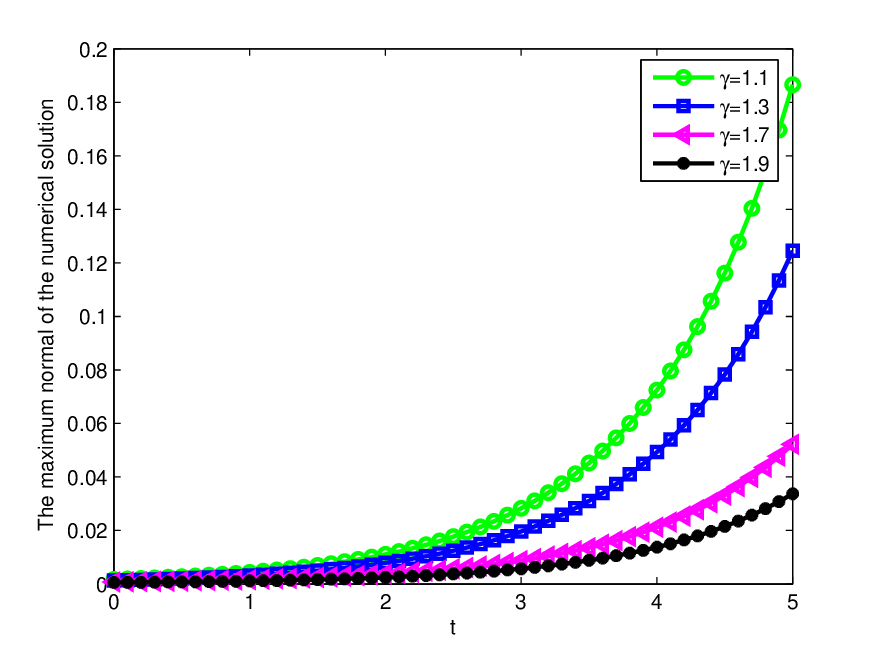}\\
  \caption{Evolutions of the maximum norm of the
  numerical solutions with different $\gamma$ values for $\tau=0.1$.}\label{Fig.33}
\end{figure}

\begin{figure}[!htb]
\centering
 \includegraphics[width=10 cm]{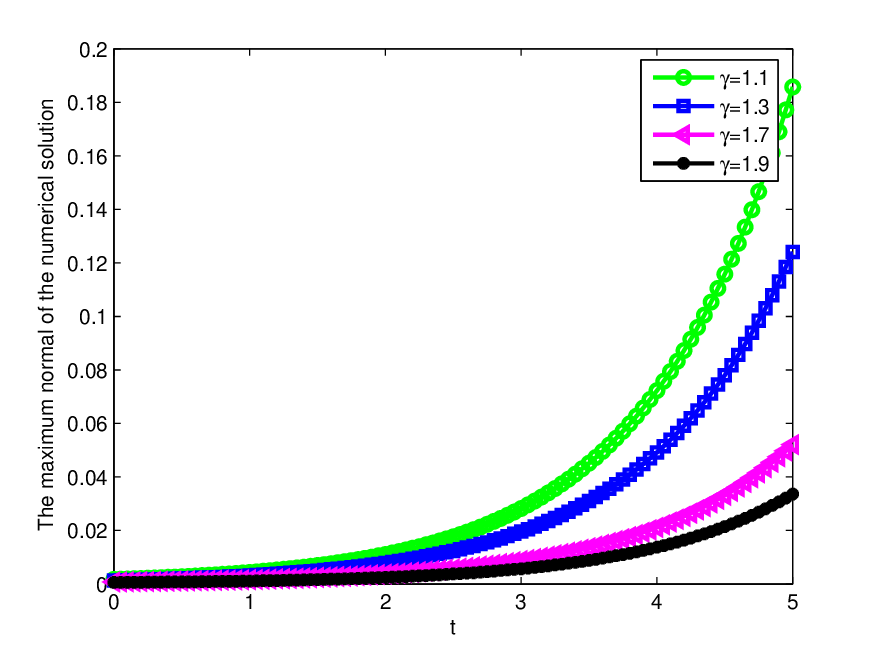}\\
  \caption{Evolutions of the maximum norm of the
  numerical solutions with different $\gamma$ values for $\tau=0.05$.}\label{Fig.34}
\end{figure}
\end{example}

\section{Conclusions}
In this article, based on the implicit integral factor method  combined with the
Pad\'{e} approximation technique, we construct an efficient numerical scheme for the
 spatial fractional Allen-Cahn equation. Subsequently, a detailed analysis
 was conducted on the unique solvability of the numerical scheme under certain
 conditions. Furthermore, we also investigated the discrete maximum principle
 preservation and energy stability under the condition $0<\tau\leq1$. The
 analysis results showed that this scheme greatly weakened the limitations of
 existing results on the selection of time step. Meanwhile, the error analysis
 under the discrete $L^\infty$ norm was also discussed in detail. Finally, a large number
 of numerical examples demonstrate the effectiveness of our algorithm
  and the correctness of theoretical analysis.

\subsection*{Declaration of competing interest}\hfill \\

The authors declare that they have no known competing financial
interests or personal relationships that could have appeared
to influence the work reported in this paper.

\subsection*{Data availability}\hfill \\

No data was used for the research described in the article.

\bibliographystyle{siamplain}

\begin{thebibliography}{99}
\bibitem{Nec2008}
{\sc Y. Nec, A. Nepomnyashchy and A. Golovin}, {\em Front-type solutions of fractional Allen-Cahn equation},
Physica D., 237 (2008), pp.~3237--3251.

\bibitem{Podlubny1999}
{\sc I. Podlubny}, {\em Fractional differential equations}. New York.: Academic Press; 1999.

\bibitem{Hou2017}
{\sc T. Hou, T. Tang and J. Yang}, {\em Numerical analysis of fully discretized
 Crank-Nicolson scheme for fractional-in-space
 Allen-Cahn equations}, J. Sci. Comput., 72 (2017), pp.~1214--1231.

\bibitem{Meerschaert2004}
{\sc M. Meerschaert, C. Tadjeran}, {\em Finite difference
approximations for fractional advection-dispersion flow
equations}, J. Comput. Appl. Math., 172 (2004), pp.~65--77.


\bibitem{Tian2015}
{\sc W. Tian, H. Zhou and W. Deng}, {\em A class of second order difference approximations
 for solving space fractional diffusion equations}, Math. Comput., 84 (2015), pp.~1703--1727.

\bibitem{Zhou2013}
{\sc H. Zhou, W. Tian and W. Deng}, {\em Quasi-compact finite
 difference schemes for space fractional diffusion equations},
  J. Sci. Comput., 56 (2013), pp.~45--66.

\bibitem{Ding2017}
{\sc H. Ding, C. Li}, {\em High-order numerical
algorithms for Riesz derivatives via constructing
new generating functions},
  J. Sci. Comput., 71 (2017), pp.~759--784.

\bibitem{Celik2012}
{\sc C. \c{C}elik, M. Duman}, {\em Crank-Nicolson method for the
fractional diffusion
equation with the Riesz fractional derivative}, J. Comput. Phys., 231 (2012), pp.~1743--1750.


\bibitem{Ortigueira2006}
{\sc M. Ortigueira}, {\em Riesz potential operators and inverses
via fractional centred derivatives}, Int. J. Math. Sci.,
 2016 (2006), pp.~1--12.

\bibitem{Ding2015}
{\sc H. Ding, C. Li and Y. Chen}, {\em High-order algorithms for Riesz derivative
 and their applications (II)}, J. Comput. Phys., 293 (2015), pp.~218--237.


\bibitem{Zhao2014}
{\sc X. Zhao, Z. Sun, Z. Hao}, {\em A fourth-order compact ADI scheme
for two-dimensional nonlinear space fractional
Schr\"{o}dinger equation}, SIAM J. Sci.
Comput., 2014 (2014). pp.~2865--2886.

\bibitem{Ding20172}
{\sc H. Ding, C. Li}, {\em High-order algorithms
 for Riesz derivative and their applications (V)}, Numer. Meth.
  Partial Diff. Equ., 33 (2017), pp.~1754--1794.


\bibitem{He2020}
{\sc D. He, K. Pan, H. Hu}, {\em A spatial fourth-order maximum principle
preserving operator splitting
scheme for the multi-dimensional fractional Allen-Cahn equation},
 Appl. Numer. Math., 151 (2020), pp.~44--63.
(2020)

\bibitem{Zhang2021}
{\sc H. Zhang, J. Yan, X. Qian, et al}, {\em On the preserving of the maximum principle and energy stability
of high-order implicit-explicit Runge-Kutta schemes for the space-fractional Allen-Cahn equation},
Numer. Algorithms. 88 (2021), pp.~1309--1336.

\bibitem{Gong2022}
{\sc Y. Gong, Y. Chen, C. Wang, et al}, {\em A new class of high-order energy-preserving schemes for the
Korteweg-de Vries equation based on the quadratic auxiliary variable (QAV) approach}, Numer. Math.
Theor. Meth. Appl., 15 (2022), pp.~ 768--792.

\bibitem{Hao2021}
{\sc Z. Hao, Z. Zhang and R. Du}, {\em Fractional centered difference scheme for
high-dimensional integral fractional
Laplacian}, J. Comput. Phys., 424 (2021), pp.~109851.

\bibitem{Xu2024}
{\sc Z. Xu, Y. Fu}, {\em Unconditional energy stability and maximum principle
preserving scheme for the Allen-Cahn equation},
Numer. Algorithms. (2024). https://doi.org/10.1007/s11075-024-01880-2.

\bibitem{Bu2019}
{\sc L. Bu, L. Mei and Y. Hou}, {\em Stable second-order schemes for the
space-fractional Cahn-Hilliard and
Allen-Cahn equations}, Comput. Math. Appl., 78 (2019), pp.~3485--3500.

\bibitem{Chen12021}
{\sc  H. Chen, H. Sun}, {\em A dimensional splitting exponential time
differencing scheme for multidimensional
fractional Allen-Cahn equations}, J. Sci. Comput., 87 (2021), pp.~1--25.

\bibitem{Chen22021}
{\sc  H. Chen, H. Sun}, {\em Second-order maximum principle preserving
Strang's splitting schemes for anisotropic
fractional Allen-Cahn equations}, Numer. Algorithms.  (2021), pp.~1--23.



\bibitem{Ding2022}
{\sc H. Ding,  Q. Yi}, {\em The construction of higher-order numerical approximation
formula for Riesz derivative and its application to nonlinear fractional differential
equations (I)}, Commun. Nonlinear Sci., 110 (2022), pp.~106394.

\bibitem{Ding2023}
{\sc H. Ding, C. Li}, {\em High-order numerical algorithm and error analysis
 for the two-dimensional nonlinear spatial fractional complex Ginzburg-Landau
 equation}, Commun. Nonlinear Sci., 120 (2023), pp.~107160.

\bibitem{Nie2008}
{\sc Q. Nie, F. Wan and Y. Zhang, et al.}, {\em Compact integration factor methods in high
spatial dimensions},  J. Comput. Phys., 2007 (2008), pp.~5238-5255.

\bibitem{Moler2003}
{\sc C. Moler, C. Loan}, {\em Nineteen dubious ways to compute the exponential
 of a matrix, twenty-five years later}, SIAM Rev., 45 (2003), pp.~3-49.

\bibitem{Chan2007}
{\sc R. Chan, X. Jin}, {\em An introduction to iterative Toeplitz solvers}, Philadelphia: SIAM; 2007.

\bibitem{Chan2011}
{\sc R. Chan}, {\em Toeplitz preconditioners for Toeplitz systems with nonnegative
 generating functions}, IMA J. Numer. Anal., 11 (1991), pp.~333-345.


\bibitem{Kincaid1991}
{\sc D. Kincaid, W. Cheney}, {\em Numerical Analysis}, Brooks/Cole Pub., California, 1991.

\bibitem{Holte2009}
{\sc J. Holte}, {\em Discrete Gr\"{o}nwall lemma and applications}, In: MAA-NCS meeting at the
 University of North Dakota, 24 (2009), pp.~1-7.

\bibitem{Akrivis1993}
{\sc  G. Akrivis}, {\em Finite difference discretization of the cubic
Schr\"{o}dinger equation}, IMA J. Numer. Anal., 13 (1993), pp.~115-124.


\end{thebibliography}

\end{document}